\documentclass[11pt]{amsart2000}
\usepackage{amscd2000,amssymb,subfigure,hyperref}
\usepackage[arrow,matrix,graph]{xy}

\numberwithin{equation}{section}

\theoremstyle{plain}
\newtheorem{theorem}[subsection]{Theorem}
\newtheorem{lemma}[subsection]{Lemma}
\newtheorem{prop}[subsection]{Proposition}
\newtheorem{corollary}[subsection]{Corollary}
\newtheorem{thm}{Theorem}

\theoremstyle{definition}
\newtheorem{remark}[subsection]{Remark}
\newtheorem{definition}[subsection]{Definition}
\newtheorem{example}[subsection]{Example}

\newtheorem*{ack}{Acknowledgments}

\newenvironment{alphenum}{

\begin{enumerate}}{\end{enumerate}}

\newenvironment{Alphenum}
{

\begin{enumerate}}{\end{enumerate}}

\newenvironment{romenum}
{

\begin{enumerate}}{\end{enumerate}}

\newenvironment{Romenum}
{

\begin{enumerate}}{\end{enumerate}}

\newcommand{\abs}[1]{\left|#1\right|}
\renewcommand{\o}[1]{\overline{#1}}

\renewcommand{\a}{\alpha}

\newcommand{\A}{\mathcal{A}}
\newcommand{\B}{\mathcal{B}}

\newcommand{\V}{\mathcal{V}}
\newcommand{\E}{\mathcal{E}}
\newcommand{\Sq}{\mathcal{S}}

\newcommand{\C}{\mathbb{C}}

\newcommand{\Z}{\mathbb{Z}}
\newcommand{\Q}{\mathbb{Q}}
\newcommand{\R}{\mathbb{R}}
\renewcommand{\P}{\mathbb{P}}

\newcommand{\wX}{\widetilde{X}}
\newcommand{\wY}{\widetilde{Y}}
\newcommand{\wA}{\widehat{\A}}

\newcommand{\OS}{\operatorname{A}^*}
\newcommand{\QOS}{\o{\operatorname{A}}^*}

\DeclareMathOperator{\rank}{rank}
\DeclareMathOperator{\coker}{coker}
\DeclareMathOperator{\id}{id}
\DeclareMathOperator{\Aut}{Aut}

\DeclareMathOperator{\ab}{ab}
\DeclareMathOperator{\codim}{codim}
\DeclareMathOperator{\pr}{pr}
\DeclareMathOperator{\im}{im}
\DeclareMathOperator{\gr}{gr}
\DeclareMathOperator{\Hilb}{Hilb}

\DeclareMathOperator{\ann}{ann}
\DeclareMathOperator{\Rad}{Rad}

\begin{document}

\title[Higher homotopy groups of arrangements]%
{Higher homotopy groups of complements of complex hyperplane
arrangements}

\author[S.~Papadima]{Stefan~Papadima}
\address{Institute of Mathematics of the Romanian Academy,
P.O. Box 1-764,
RO-70700 Bucharest, Romania}
\email{\href{mailto:Stefan.Papadima@imar.ro}{%
Stefan.Papadima@imar.ro}}
\urladdr{\href{http://www.imar.ro/~spapadim}{%
http://www.imar.ro/\~{}spapadim}}

\author[A.~I.~Suciu]{Alexander~I.~Suciu}
\address{Department of Mathematics,
Northeastern University,
Boston, MA 02115, USA}
\email{\href{mailto:alexsuciu@neu.edu}{%
alexsuciu@neu.edu}}
\urladdr{\href{http://www.math.neu.edu/~suciu/}{%
http://www.math.neu.edu/\~{}suciu}}




\begin{abstract}
We generalize results of Hattori on the topology of complements
of hyperplane arrangements, from the class of generic arrangements,
to the much broader class of hypersolvable arrangements.
We show that the higher homotopy groups of the complement
vanish in a certain combinatorially determined range, and we
give an explicit $\Z\pi_1$-module presentation of $\pi_p$,
the first non-vanishing higher homotopy group. We also
give a combinatorial formula for the $\pi_1$-coinvariants
of $\pi_p$.  

For affine line arrangements whose cones are hypersolvable, 
we provide a minimal resolution of $\pi_2$, and study some 
of the properties of this module.  For graphic arrangements 
associated to graphs with no $3$-cycles, we obtain information 
on $\pi_2$, directly from the graph.  The $\pi_1$-coinvariants 
of $\pi_2$ may distinguish the homotopy $2$-types of arrangement 
complements with the same $\pi_1$, and the same Betti numbers 
in low degrees.
\end{abstract}

\maketitle

\tableofcontents

\section{Introduction}
\label{sec:intro}

\subsection{Background}
\label{subsec:back}
One of the fundamental problems in the topological study of
polynomial functions, $f\colon (\C^{\ell},{\mathbf{0}})\to (\C,0)$,
is the computation of the homotopy groups of the complement
to the hypersurface $V(f)=f^{-1}(0)$. A well-known algorithm for
finding
a finite presentation for $\pi_1(\C^{\ell}\setminus V(f))$
was given by Zariski and VanKampen in the early 1930's. Much less
is known about the higher homotopy groups of the complement,
except when $V(f)$ is irreducible, in which case
the Zariski-VanKampen method can be extended to
give information about $\pi_k(\C^{\ell}\setminus V(f))\otimes\Q$,
for $k>1$, see~\cite{Li2}.

In this paper, we concentrate on the simplest kind of polynomial $f$
for which the hypersurface $V(f)$ is {\em not} irreducible.
Namely, suppose $f$ factors completely into distinct, degree~one
factors.
Then $f$ is the defining polynomial of a hyperplane arrangement, $\A$,
with union $V(f)=\bigcup_{H\in \A} H$, and complement
$X(\A)=\C^{\ell}\setminus \bigcup_{H\in \A} H$. The cohomology ring
of $X=X(\A)$ was computed by Brieskorn \cite{Br}.
Orlik and Solomon \cite{OS} expressed $H^*(X)$ in terms of the
combinatorics of $\A$, encoded in the intersection lattice, $L(\A)$.
In particular, the Poincar\'{e} polynomial,
$P_{\A}(T) =\sum_{k=1}^{\ell} b_k(X) T^k$,
admits a simple combinatorial expression, see Orlik and
Terao~\cite{OT}.
On the other hand, the fundamental group of the complement, $\pi_1(X)$,
is not determined by $L(\A)$ alone, as the example of
Rybnikov~\cite{Ryb} shows.

For certain arrangements, all the higher homotopy groups of the
complement vanish. Examples of such $K(\pi,1)$ arrangements include
the {\em simplicial} arrangements (Deligne~\cite{De}),
and the {\em supersolvable} arrangements (Terao~\cite{Te}).
Examples of non-$K(\pi,1)$ arrangements,
and methods for detecting the first non-vanishing higher
homotopy group of their complements,
were given by Falk~\cite{Fa2} and Randell~\cite{Ra2}
(see also the recent survey \cite{FR3}).

The first (and, up to now, only) explicit computation of
non-trivial higher homotopy groups of arrangement complements
was made by Hattori \cite{Ha}.
An arrangement $\A$ in $\C^{\ell}$, $\ell>1$,
is called {\em generic} if, for all $\B\subset \A$, the
intersection $\bigcap_{H\in \B} H$ has codimension $\abs{\B}$
when $\abs{\B}\le \ell$, and is empty when $\abs{\B}> \ell$.
The standard example is the {\em Boolean} arrangement
of coordinate hyperplanes in $\C^{n}$,
with complement $(\C^*)^{n}$.
Hattori used the minimal cell decomposition of
$(\C^*)^{n} \simeq T^{n}$ to find
an explicit, minimal cell decomposition for the complement
of an arbitrary generic arrangement. More precisely, if $\A$
is an arrangement of $n$ hyperplanes in general position in $\C^{\ell}$
($n>\ell$), then $X(\A)\simeq (T^{n})^{(\ell)}$.
  From this decomposition, Hattori deduced:
\begin{Alphenum}
\item\label{hat1} $\pi_1(X)=\Z^n$.
\item\label{hat2} $\pi_k(X)=0$ for $1< k < \ell$.
\item\label{hat3} $\pi_{\ell}(X)$ admits a free
$\Z\pi_1$-resolution of length $n-\ell$.
\end{Alphenum}

The simplest example is that of $3$ generic affine lines in $\C^2$.
In that case, the complement $X$ has the homotopy type of
the $2$-skeleton of the $3$-torus,
$T^3=K(\Z^3, 1)$. Looking at the universal cover, $\wX$, we thus see
that
$\pi_2(X)=H_2(\wX)$ is a free $\Z\Z^3$-module, generated by the
boundary of
a cubical $3$-cell from the standard decomposition of
$\widetilde{T^3}=\R^3$.

\subsection{Results}
\label{subsec:basic}
In this paper, we set out to generalize Hattori's results to the wider
class of {\em hypersolvable arrangements}. This combinatorially defined
class, introduced in
\cite{JP1}, includes both supersolvable arrangements and (cones of)
generic arrangements.

A hypersolvable arrangement $\A$ admits a ``supersolvable deformation,"
$\wA$, which preserves the collinearity relations.
For example, if we start with $n\ge 3$ generic lines in $\C^2$, and
take $\A$ to be the respective central arrangement of planes in $\C^3$,
then $\wA$ is the Boolean arrangement in $\C^n$.
In general, $X(\A)$ has the same fundamental
group as $X(\wA)$; see \cite{JP1, JP2}.
Moreover, $\pi_1(X(\wA))$ is a (special kind of) iterated
semidirect product of finitely-generated free groups;
see \cite{Te, FR, CSbm} and Theorem~\ref{thm:hypsolv}.
These facts together provide the generalization of
Hattori's result \eqref{hat1} to hypersolvable arrangements.

The key tool for generalizing
\eqref{hat2} and \eqref{hat3} to complements of
  hypersolvable arrangements is the existence of {\em minimal} cell
structures, on both $X(\A)$ and $X(\wA)=K(\pi_1(X(\A)), 1)$.

To find a presentation for the first higher non-vanishing
homotopy group, we thus turn to a general study of minimal cell
decompositions.
The idea is to get higher homotopy information on
a connected, finite-type, CW-space $X$, by
comparing it to its classifying space $K(\pi,1)$,
where $\pi =\pi_1(X)$.
We are thus led to introduce a homotopy-type invariant of $X$,
called the {\em order of $\pi_1$-connectivity}, which measures
the rational-homology deviation of $X$ from asphericity:
\[
p(X):=\sup\{q \mid b_r(X)=b_r(K(\pi,1)),\ \forall r\le q\}.
\]
(If $X$ is $1$-connected, and $H_*(X)$ is torsion-free, then $p(X)$ is
the
usual order of connectivity of $X$.)

A (connected, finite-type) CW-space $X$ is said to be {\em minimal} if
it has a CW-structure with $b_k$ $k$-cells, for all $k$, where $b_k$ is
the $k$-th Betti number of $X$.
In \S\ref{sec:general}, we will prove the following result.

\begin{theorem}
\label{thm:intro2}
Let $X$ and $Y$ be two minimal CW-complexes, with
cohomology rings
generated in degree $1$. Set $\pi=\pi_1(X)$ and $p=p(X)$.
Assume that $Y$ is aspherical and that there is a cellular map,
$j\colon X\to Y$, such that $j\!\left\vert_{X^{(p)}}\right.=\id$.
Then $j$ is a classifying map. Moreover:
\begin{enumerate}
\item \label{z1} $X$ is aspherical if and only if $p=\infty$.
\item \label{z2} If $p< \infty$, then the first non-trivial
higher homotopy group of $X$ is $\pi_p(X)$, which
has the following finite
$\Z\pi$-presentation:
\begin{equation}
\label{eq:intropres}
(H_{p+2}(Y) \oplus H_{p+1}(X))\otimes \Z\pi
\xrightarrow{D_p := \partial_{p+2}+ \tilde{j}_{p+1}}
H_{p+1}(Y)\otimes\Z\pi\to \pi_p(X)\to 0,
\end{equation}
where $\partial_*$ is the differential of the $\pi$-equivariant
cellular chain complex of the universal cover of $Y$, and
$\tilde{j}_*$ is the $\Z\pi$-chain map induced by the lift of $j$
to universal covers.
\item \label{z3} If $p< \infty$, then the group of $\pi$-coinvariants
of $\pi_p(X)$ is isomorphic to $H_{p+1}(Y,X)$, where
$H_*(Y,X) := \coker(j_* \colon H_*(X) \to H_*(Y))$.
\end{enumerate}
\end{theorem}

The theorem provides a complete generalization of Hattori's
result \eqref{hat2} in this setting, and a partial generalization
of \eqref{hat3}. If $\dim X= p(X)$ and $Y$ is a finite complex,
the presentation \eqref{eq:intropres} extends to a finite-length, free
$\Z\pi$-resolution of $\pi_p(X)$, see Remark~\ref{rem:hatres}.
In particular, if $X$ is the complement of a generic arrangement,
then $p(X)=\ell$, and
\eqref{eq:intropres} may be continued to Hattori's
resolution \eqref{hat3}.

In \S\ref{sec:fitting}, we follow a standard approach and
extract from the above presentation
of $\pi_p(X)$ more manageable invariants of homotopy type:
the subvarieties of the complex torus $(\C^*)^n$,
$n=b_1(X)$, defined by the Fitting ideals of
$\pi_p(X)\otimes_{\Z\pi}\Z\Z^n$. In turn, we identify these
varieties with the jumping loci for homology with coefficients in
rank~$1$ local systems of the pair $(K(\pi,1),X)$.

We now return to the case where $X=X(\A)$ is the complement of an
(essential) arrangement $\A$. From recent work of Dimca \cite{D} and
Randell \cite{R}, we know that $X$ is minimal.
We also know (from \cite{Br}) that $H^*(X)$ is generated in degree $1$.
Since $X$ may fail to possess any finite-type $K(\pi,1)$,
our approach does not work in this generality.
If $\A$ is hypersolvable, though,
we may take $K(\pi,1)=X(\wA)$, where $\wA$ is the
supersolvable deformation of $\A$.

We devote \S\ref{sec:hypsolv} to the description of
the first non-vanishing higher homotopy group of a hypersolvable
(non-supersolvable) arrangement, $\A$. The combinatorics and
the homotopy theory of the complement of $\A$ do not change, when
passing to the associated essential arrangement, $\A_{\rm ess}$;
see \cite{OT}. Let $\wA$ be the supersolvable deformation of the
essential hypersolvable arrangement $\A_{\rm ess}$.
Exploiting  the method from \cite[Section 5]{DP},
we show in Theorems \ref{thm:minarr} and \ref{thm:hypsolv} 
\eqref{part0}
how to replace, up to homotopy, the inclusion
$J\colon X(\A_{\rm ess}) \hookrightarrow X(\wA)$ by a map,
$j\colon X \to Y$, which satisfies {\em all} hypotheses from
Theorem~\ref{thm:intro2}. This leads to the following result.

\begin{theorem}
\label{thm:intro1}
Let $\A$ be a hypersolvable arrangement, with complement $X=X(\A)$,
fundamental group $\pi=\pi_1(X)$, and order of $\pi_1$-connectivity
$p=p(X)$. Then:
\begin{enumerate}
\item \label{t1} $X$ is aspherical $\Longleftrightarrow$ $\A$ is
supersolvable
$\Longleftrightarrow$ $p= \infty$.
\item \label{t2} If $p< \infty$, then the first
non-vanishing higher homotopy group of $X$
is $\pi_p(X)$, with finite $\Z\pi$-presentation
\eqref{eq:intropres}.
\item \label{t3} If $p< \infty$, then the group of
$\pi$-coinvariants of $\pi_p(X)$ is free
abelian, with (strictly positive) combinatorially determined rank.
\end{enumerate}
\end{theorem}

In Corollary~\ref{cor:pcomb}\eqref{i1}, we give a combinatorial
interpretation for $p(X)$.
The precise formula for the coinvariants is provided by
Theorem~\ref{thm:hspipres}\eqref{p3}. A similar formula was obtained by
Randell \cite{Ra2}, for generic sections of aspherical arrangements---%
a class of arrangements which overlaps with the hypersolvable class,
but neither includes it, nor is included in it.

\subsection{Applications}
\label{subsec:appl}
Particularly simple is the case of affine line arrangements
in $\C^2$. These arrangements represent both the simplest
case of non-irreducible plane algebraic curves, and the simplest
case of hyperplane arrangements. As such, they have been the object
of intense investigation, see e.g. \cite{Fa2, CSbm, Fn, Ryb, FR3}.
For one, the fundamental group of an arbitrary hyperplane
arrangement complement can be identified with the fundamental group
of an affine line arrangement (by the Hamm-L\^{e} theorem), thereby
making line arrangements key to the understanding of all
arrangements. For another, complements of affine line arrangements
need not be aspherical (unlike, say, complements of
weighted-homogeneous plane curves, which always are), thereby
making for a richer object of topological study.

In \S\ref{sec:pi2pres}, we consider affine line arrangements
whose cones are hypersolvable. In Theorem~\ref{thm:pi2hyper},
we go further, providing a minimal, finite-length resolution
for $\pi_2$, which completely generalizes Hattori's resolution
\eqref{hat3} in this context. We also obtain some finer
information about the $\Z\pi_1$-module structure of $\pi_2$:
it is neither projective (except in a very special,
combinatorially decidable case, when it is free, with
rank combinatorially determined), nor nilpotent (except
if it is trivial).

Another class of arrangements which can be fairly well understood
from our point of view is that of (hypersolvable) graphic
arrangements. In \S\ref{sec:graphics}, we implement in
this setting our method for higher homotopy computations.
The class of arrangements associated
to graphs without $3$-cycles provides a natural, rich supply
of hypersolvable arrangements, which are neither supersolvable
nor generic, and for which homotopy information may be extracted
directly from the graph. As an illustration, we exhibit two
graphic arrangements, whose complements have the same $\pi_1$,
but non-isomorphic $\pi_2$'s (when viewed as $\Z\pi_1$-modules).

This version of our paper corrects an oversight from
\cite{PS}.

\begin{ack}
This research was supported by the Volkswagen-Stiftung
(RiP-pro\-gram at Oberwolfach). The authors wish to thank
the Mathematisches Forschungs\-institut for the warm
hospitality and excellent facilities provided.
The second author was also supported by an
RSDF grant from Northeastern University.
\end{ack}

\section{Minimal cell decompositions and homotopy groups}
\label{sec:general}

\subsection{Minimal cell decompositions}
\label{subsec:min}
Given a space $X$, consider the following conditions on its homotopy
type:
\begin{romenum}
\item \label{cond:finite}
$X$ is homotopy equivalent to a connected, finite-type CW-complex;
\item \label{cond:free}
The integral homology groups $H_*(X)$ are torsion-free;
\item \label{cond:surj}
The cup-product map $\cup\colon\bigwedge^* H^1(X)\to H^*(X)$ is 
surjective.
\end{romenum}

These three conditions abstract some well-known topological properties
of complements of complex hyperplane arrangements. Next, we delineate
a class of spaces that satisfy condition \eqref{cond:finite} and
a much stronger form of \eqref{cond:free}.

\begin{definition}
\label{def:min}
A space $X$ is called {\em minimal} if $X$ has the homotopy type
of a connected, finite-type, CW-complex $K$ such that
\begin{equation}
\label{eq:minimal}
\#\{\text{$q$-cells in $K$}\}= \rank H_q(X ; \Z),\ \text{for all $q\ge
0$}.
\end{equation}
\end{definition}

This definition implies at once that all the (abelian) groups $H_q(X)$
are finitely-generated and torsion-free. Consequently, we may
unambiguously speak about the Betti numbers of $X$, $b_q(X)$,
without specifying the coefficients.

Let $X$ be a minimal space, and let
$C_*(X)$ be the cellular chain complex of $X$,
corresponding to a minimal CW-decomposition.
Let $\pi=\pi_1(X)$ be the fundamental group, $\Z\pi$ its group ring,
and $\epsilon\colon\Z\pi\to\Z$ the augmentation map. Let $\wX$ be the
universal cover of $X$, and let $(C_*(\wX),d_*)$ be the
$\pi$-equivariant
chain complex of $\wX$, with $C_q(\wX)=C_q(X)\otimes\Z\pi$
and $d_q\colon C_q(\wX)\to C_{q-1}(\wX)$. (Note that, when turning left
$\Z\pi$-modules into right $\Z\pi$-modules, one has to
replace the action of $x\in \pi$ by that of $x^{-1}$.)
By minimality,
all the boundary maps $d_q$ are {\em $\epsilon$-minimal},
in the sense that $d_q\otimes_{\Z\pi}\Z=0$.

\begin{example}
\label{ex:tori}
The standard example of a space admitting a minimal cell decomposition
is the $n$-torus, $T^n$. Identifying $\pi_1(T^n)=\Z^n$, with basis
$\{x_i\}_{i}$, and $C_q(T^n)=\bigwedge^q \Z^n$, with basis
$\{\sigma_I=\sigma_{i_1}\cdots \sigma_{i_q}\}_{I}$, the boundary map
$d_q\colon \bigwedge^{q}\Z^n\otimes \Z\Z^n\to\bigwedge^{q-1} 
\Z^n\otimes
\Z\Z^n$ is given by $d_{q}(\sigma_I)=\sum_{r=1}^{q} (-1)^{r-1}
\sigma_{I\setminus\{i_r\}}\otimes (x^{-1}_{i_{r}}-1)$.
\end{example}

\begin{example}
\label{ex:semiprod}
More generally, let $\pi=F_{d_n}\rtimes_{\rho_{n-1}} F_{d_{n-1}}\rtimes
\cdots \rtimes_{\rho_1} F_{d_1}$ be an iterated semidirect product of 
free
groups, with $\rho_i$ acting as the identity in homology, and 
$X=K(\pi,1)$
a corresponding Eilenberg-MacLane space.  A finite, minimal cell
decomposition of $X$ is given in \cite{CScc}:
The number of cells is read off the Poincar\'e polynomial,
$P_X(T)=\prod_{i=1}^{n} (1+d_i T)$, and the ($\epsilon$-minimal)
boundary maps, $d_q\colon C_q(\wX)\to C_{q-1}(\wX)$, are given 
explicitly
in terms of Fox Jacobians of the monodromy operators $\rho_i$,
see \cite[ Thm.~2.10, Prop.~3.3, and Cor.~3.4]{CScc}.
\end{example}

\begin{remark}
\label{rem:knotslinks}
Not all manifolds admit minimal cell decompositions.
For example, if $X$ is the complement of a non-trivial
knot in $S^3$, then $X$ has no minimal cell decomposition,
not even up to $q=1$. See also the monograph by
Sharko~\cite{Sh} for various other definitions of minimality
in related contexts.
\end{remark}

Now assume $X$ is a minimal space for which there exists a
minimal Eilenberg-MacLane space $Y=K(\pi,1)$.
Let $j\colon X\to Y$ be a classifying map.
Without loss of generality, we may assume $j$ respects the given
(minimal) CW-decompositions on $X$ and $Y$. Then the chain
map $j_*\colon C_*(X)\to C_*(Y)$ lifts to an equivariant chain map
$\tilde{j}_*\colon (C_*(\wX),d_*) \to (C_*(\wY), \partial_*)$,
\[
\tilde{j}_* = \{ \tilde{j}_{q}\colon  H_{q}(X)\otimes\Z\pi\rightarrow
H_{q}(Y)\otimes\Z\pi \}_{q\ge 0} \, .
\]

\subsection{Homotopy groups}
\label{subsec:homotopy}
We now analyze the homotopy groups of certain minimal spaces.
In order to state our results,
we need to introduce one more notion.

\begin{definition}
\label{def:px}
Let $X$ be a space satisfying condition~\eqref{cond:finite}. Define
the {\em order of $\pi_1$-connec\-tivity} of $X$ to be
\[
p(X):=\sup\{q \mid b_r(X; \Q)=b_r(K(\pi_1(X),1); \Q),\ \forall r\le 
q\}.
\]
\end{definition}

\begin{remark}
\label{rem:px}
The terminology is borrowed from the simply-connected case:
if $\pi_1(X)=0$ and $X$ also satisfies \eqref{cond:free}, then
$p(X)$ is the usual order of connectivity of $X$.
Note that $p(X)$ is a positive integer, depending only on the
homotopy type of $X$. Furthermore, if $X$
satisfies conditions \eqref{cond:finite}--\eqref{cond:surj},
then $p(X)\ge 2$.
\end{remark}

\begin{remark}
\label{rem:splitinj}
Set $Y=K(\pi_1(X),1)$, and consider a classifying map,
$j\colon X\to Y$.  Assume both $X$ and $Y$ satisfy conditions
\eqref{cond:finite}--\eqref{cond:surj} from~\ref{subsec:min}.
These conditions readily imply that $j$
induces a split surjection on cohomology, and a split injection
on homology. Consequently, $j_{r}\colon H_r(X)\to H_r(Y)$ is an
isomorphism, for all $r\le p(X)$, and the groups
$H_*(Y,X):=\coker(j_{*}\colon  H_*(X)\to H_*(Y))$
fit into split exact sequences
\begin{equation}
\label{eq:jpi}
0\to H_*(X)\xrightarrow{j_*}H_*(Y)\xrightarrow{\Pi_*}H_*(Y,X)\to 0.
\end{equation}
\end{remark}

The next two results provide a complete proof of
Theorem~\ref{thm:intro2} from the Introduction.
\begin{theorem}
\label{thm:hure}
Let $X$ and $Y$ be two minimal CW-complexes satisfying
condition \eqref{cond:surj} from
\S\ref{subsec:min}. Set $\pi =\pi_1(X)$ and $p=p(X)$.
Assume that $Y$ is aspherical, and that there is a cellular map,
$j\colon X\to Y$, such that the restriction of $j$ to the $p$-skeleton,
$X^{(p)}$, is the identity. Then $j$ is a classifying map, and:
\begin{enumerate}
\item \label{f1} $\wX$ is $(p-1)$-connected.
\item \label{f2} If $p< \infty$, then
\begin{equation}
\label{eq:pippres}
(H_{p+2}(Y)\otimes \Z\pi) \oplus (H_{p+1}(X)\otimes \Z\pi)
\xrightarrow{D_p := \partial_{p+2}+ \tilde{j}_{p+1}}
H_{p+1}(Y)\otimes\Z\pi\to \pi_p(X)\to 0
\end{equation}
is a finite presentation of $\pi_p(X)$ as $\Z\pi$-module.
\end{enumerate}
\end{theorem}
\begin{proof}
Since $p\ge 2$, $Y$ is a $K(\pi, 1)$ and $j$ is a classifying map.

\eqref{f1} We have to show that $\pi_q(\wX)=0$, for $q< p$.
Of course, $\pi_1(\wX)=0$. Fix $1<q<p$, and assume that
$\pi_r(\wX)=0$, for $r<q$. By the Hurewicz isomorphism theorem,
$\pi_q(\wX)=H_q(\wX)$. By minimality of $X$ and $Y$,
we have a commuting ladder between the (equivariant)
chain complexes of $\wX$ and $\wY$:
\begin{equation*}
\label{eq:ladder1}
\begin{CD}
C_*(\wX):\qquad &
H_{q+1}(X)\otimes\Z\pi @>d_{q+1}>>
H_{q}(X)\otimes\Z\pi @>d_{q}>>
H_{q-1}(X)\otimes\Z\pi \\
&
@VV{\tilde{j}_{q+1}}V
@VV{\tilde{j}_{q}}V
@VV{\tilde{j}_{q-1}}V \\
C_*(\wY):\qquad &
H_{q+1}(Y)\otimes\Z\pi @>\partial_{q+1}>>
H_{q}(Y)\otimes\Z\pi @>\partial_{q}>>
H_{q-1}(Y)\otimes\Z\pi
\end{CD}
\end{equation*}
The three vertical arrows are isomorphisms,
since $\tilde{j}_{r}= \id$, for $r\le p(X)$.
It follows that $H_q(\wX)=H_q(\wY)$. But $\wY$ is acyclic,
and so $\pi_q(\wX)=0$.

\eqref{f2} Consider the commuting diagram
\begin{equation*}
\label{eq:ladder2}
\begin{CD}
&&
H_{p+1}(X)\otimes\Z\pi @>d_{p+1}>>
H_{p}(X)\otimes\Z\pi @>d_{p}>>
H_{p-1}(X)\otimes\Z\pi\\
&&
@VV{\tilde{j}_{p+1}}V
@VV{\id}V
@VV{\id}V \\
H_{p+2}(Y)\otimes\Z\pi @>\partial_{p+2}>>
H_{p+1}(Y)\otimes\Z\pi @>\partial_{p+1}>>
H_{p}(Y)\otimes\Z\pi @>\partial_{p}>>
H_{p-1}(Y)\otimes\Z\pi\\
\end{CD}
\end{equation*}
By Part~\eqref{f1} and Hurewicz, $\pi_p(X)=H_p(\wX)$.
A diagram chase yields isomorphisms
\[
H_p(\wX)=\frac{\ker d_p}{\im d_{p+1}}=
\frac{\ker \partial_p}{\im d_{p+1}}=
\frac{\im\partial_{p+1}}{\im
(\partial_{p+1}\circ \tilde{j}_{p+1})}
\xleftarrow{\partial_{p+1}}
\frac{H_{p+1}(Y)\otimes \Z\pi}{\im \partial_{p+2} + \im
\tilde{j}_{p+1}} =
\coker (D_p)\, ,
\]
and we are done.
\end{proof}

\begin{corollary}
\label{cor:coinv}
With assumptions as above,
and if $p=p(X)<\infty$, then the group of coinvariants
of $\pi_p(X)$ under the action of $\pi=\pi_1(X)$ is given by
\[
(\pi_p(X))_{\pi}=H_{p+1}(Y,X).
\]
In particular, $(\pi_p(X))_{\pi}\ne 0$.
\end{corollary}
\begin{proof}
 From presentation \eqref{eq:pippres}, we compute
$(\pi_p(X))_{\pi}= \coker (D_p \otimes_{\Z\pi} \Z)$.
At the same time, $\coker (D_p \otimes_{\Z\pi} \Z)=\coker (\tilde{j}_{p+1} \otimes_{\Z\pi} \Z)$,
since $\partial_{p+2}$ is $\epsilon$-minimal, and $\tilde{j}_{p+1}
\otimes_{\Z\pi} \Z = j_{p+1}$, since $\tilde{j}_*$ lifts $j_*$. This
proves the first claim. To finish the proof, note that
$H_{p+1}(Y, X)\ne 0$, by the definition of $p$.
\end{proof}

Assuming $p< \infty$ in Theorem~\ref{thm:hure}, one may
also consider the following $\Z \pi$-linear map:
\begin{equation}
\label{eq:deltamin}
H_{p+2}(Y)\otimes \Z\pi
\xrightarrow{\Delta_{p} := (\Pi_{p+1}\otimes \id)\circ \partial_{p+2}}
H_{p+1}(Y, X)\otimes \Z\pi \, .
\end{equation}
In the next two remarks, we present some cases when the presentation 
matrix
$D_{p}$ from \eqref{eq:pippres} may be replaced by the
(simpler, $\epsilon$-minimal) presentation matrix $\Delta_{p}$ from
\eqref{eq:deltamin}.

\begin{remark}
\label{rem:hatres}
If $X$ has the homotopy type of a CW-complex of dimension
$p=p(X)$, and $Y$ is finite, then the presentation \eqref{eq:pippres}
for
$\pi_p(X)$ may be continued to a free $\Z\pi$-resolution of
length $d-p$, where $d=\dim Y$. (We shall encounter such a
situation later on, in Theorem~\ref{thm:pi2hyper}.)
Indeed, let $(C_*(\wY),\partial_*)$ be the $\pi$-equivariant
cellular chain complex of $\wY$.
Note that $H_{p+1}(X)=0$ (since $\dim X= p$), and so
$H_{p+1}(Y,X)=H_{p+1}(Y)$ and $\Delta_{p}=\partial_{p+2}=D_{p}$.
Hence, $\pi_p(X)$ has finite, free, $\epsilon$-minimal, resolution
\begin{equation*}
\label{eq:pipxres}
0\to C_{d}(\wY)\xrightarrow{\partial_{d}}
C_{d-1}(\wY)\to\cdots\to
C_{p+2}(\wY)
\xrightarrow{\partial_{p+2}}
C_{p+1}(\wY)\to \pi_p(X)\to 0.
\end{equation*}
\end{remark}

\begin{remark}
\label{rem:subcomplex}
An especially simple situation where Theorem~\ref{thm:hure}
applies is as follows.  Let $Y$ be a minimal $K(\pi,1)$-complex
satisfying condition \eqref{cond:surj},
and let $X\subset Y$ be a proper, connected subcomplex,
such that $X^{(2)}=Y^{(2)}$.
Since $Y$ is minimal, $X$ is also minimal, and \eqref{cond:surj} also
holds for $X$. Since $X$ and $Y$ share the
same $2$-skeleton, the inclusion $j\colon X\to Y$ is a classifying map,
inducing an isomorphism on $\pi_1$.
Since $X$ is a subcomplex of $Y$, we have an exact sequence
of cellular chain complexes,
$0\to C_*(X)\xrightarrow{j_*} C_*(Y) \xrightarrow{\pr_*} C_*(Y,X)\to 
0$,
and $\Pi_*=\pr_*\colon H_*(Y)\to H_*(Y,X)$. Moreover,
\[
p(X)=\max\{ q \mid
\#\{\text{$r$-cells of $X$}\}=
\#\{\text{$r$-cells of $Y$}\}, \
\forall r\le q\},
\]
and $p(X) <\infty$. Therefore, $\pi_k(\wX)=0$, for $k<p=p(X)$,
and $\pi_p(X)$ has a finite $\Z\pi$-presentation,
given in \eqref{eq:pippres}. Obviously,
$\tilde{j}_{p+1}= j_{p+1}\otimes \id$, since $X^{(p+1)}$ is
a subcomplex of $Y^{(p+1)}$. Hence,
$\Pi_{p+1}\otimes \id$ induces an isomorphism between
$\coker (D_{p})$ and $\coker (\Delta_{p})$.
\end{remark}

\begin{remark}
\label{rem:1conn}
Most of the results in this section have only relevance for
non-simply-con\-nected spaces. Indeed, if $X$ has the
homotopy type of a finite-type CW-complex, and
$\pi_1(X)=0$, then $X$ cannot satisfy condition~\eqref{cond:surj},
unless $X$ is contractible.

On the other hand, if $X$ is $1$-connected and satisfies conditions
\eqref{cond:finite} and \eqref{cond:free}, then $X$ has the homotopy
type of a minimal CW-complex $K$. The complex $K$ may be obtained
from a bouquet of spheres $\bigvee_{b_{p+1}(X)} S^{p+1}$, where
$p=p(X)$, by attaching suitable cells. For details of the proof, see
Anick~\cite{An}, where the notion of minimality for simply-connected
spaces was first introduced.
\end{remark}

\section{Fitting ideals and jumping loci}
\label{sec:fitting}

The $\Z\pi$-module $\pi_p(X)$ determined in Theorem~\ref{thm:hure}
can be rather intractable. We now associate to $\pi_p(X)$
a more manageable module (over a commutative ring), and extract from
it invariants that can be understood as jumping loci for
homology with coefficients in rank~$1$ local systems.

\subsection{Fitting ideals}
\label{subsec:fitt}
Let $\pi$ be a group, with abelianization $\pi^{\ab}\cong \Z^n$,
and let $M$ be a finitely-presented module over $\Z\pi$.
Let $\widetilde{M}=M\otimes_{\Z\pi} \Z\Z^n$
be the module over $\Z\Z^n$ obtained by extending
scalars via $\Z\pi\xrightarrow{\ab}\Z\Z^n$.
For $k\ge 0$, let $F_k(\widetilde{M})$ be the corresponding
{\em $k^{\text{th}}$ Fitting ideal},
generated by the codimension $k-1$ minors of a presentation
matrix for $\widetilde{M}$. As is well-known, the Fitting ideals
are independent of the choice of presentation, see
e.g.~\cite[p.~493]{Ei}.
Now fix a basis $\{x_1,\dots, x_n\}$ for $\Z^n$. Then, the group
ring $\C\Z^n=\Z\Z^n\otimes \C$ may be identified with
$\C[x_1^{\pm 1},\dots ,x_n^{\pm 1}]$,
the coordinate ring of the complex algebraic torus $(\C^*)^n$.
For each $k\ge 0$, the $k^{\text{th}}$ Fitting ideal of
$\widetilde{M}$ defines a subvariety of this torus,
\[
V_k(M):=\{t \in (\C^*)^n \mid g(t)=0, \
\forall g\in F_k(\widetilde{M})\otimes \C\}.
\]
Alternatively, $V_k(M)$ can be described as the
variety defined by the annihilator of $\bigwedge^k \widetilde{M}$.
Indeed,
$\Rad(F_k(\widetilde{M}))=\Rad(\ann(\bigwedge^k \widetilde{M}))$,
see \cite[pp.~511-513]{Ei}.

Let $X$ be a path-connected space. Assume that $\pi=\pi_1(X)$
has abelianization $\Z^n$, and that $\pi_p(X)$ is a finitely
presented $\Z\pi$-module.

\begin{definition}
\label{def:vpi}
The {\em $k$-th Fitting variety} of $\pi_p(X)$ is the subvariety
of the complex algebraic $n$-torus, $V_k(\pi_p(X))$, defined as above.
\end{definition}

A standard argument (which we include for the sake of completeness)
shows that the Fitting varieties are invariants of the
$\pi_1$-module $\pi_p$.

\begin{prop}
\label{prop:pipinv}
Suppose $f\colon \pi_1(X)\to \pi_1(X')$ and $g\colon \pi_p(X)\to 
\pi_p(X')$
are compatible isomorphisms, i.e., $g(x m)=f(x) g(m)$,
for all $x\in \pi_1(X)$ and $m\in \pi_p(X)$. Then there
is a monomial isomorphism $\Phi\colon (\C^*)^n\to (\C^*)^n$ such that
$\Phi(V_k(\pi_p(X))=V_k(\pi_p(X'))$.
\end{prop}

\begin{proof}
The extension of $f_{\ab}\colon H_1(X) \to H_1(X')$ to group rings
gives rise to an isomorphism between
$\widetilde{M}=\pi_p(X)\otimes_{\Z\pi_1(X)} \Z H_1(X)$ and
$\widetilde{M'}=\pi_p(X')\otimes_{\Z\pi_1(X')} \Z H_1(X')$,
and thus maps
$\Rad(\ann(\bigwedge^k \widetilde{M}))$ bijectively to
$\Rad(\ann(\bigwedge^k \widetilde{M'}))$.

Now identify $H_1(X)$ and $H_1(X')$ with $\Z^n$, and let
$\phi=\left(\phi_{i,j}\right)\colon \Z^n\to\Z^n$
be the matrix of $f_{\ab}$ under this identification.
Let $\Phi\colon (\C^*)^n\to (\C^*)^n$ be the
corresponding monomial isomorphism, given by $\Phi(t_i)=
t_1^{\phi_{i,1}} \cdots t_n^{\phi_{i,n}}$.
It is readily verified that $\Phi$ preserves the
Fitting varieties.
\end{proof}

\begin{corollary}
\label{cor:fittinv}
For each $k\ge 0$, the monomial isomorphism type of $V_k(\pi_p(X))$ is
a
homotopy type invariant for $X$. \hfill\qed
\end{corollary}

Now let $X$ be a space satisfying the conditions of
Theorem~\ref{thm:hure}\eqref{f2}.
Then $\pi^{\ab}\cong \Z^n$, where $n=b_1(X)$. Let $p=p(X)$ be the
order
of $\pi_1$-connectivity of $X$. The Fitting ideals of
$\pi_p(X)\otimes_{\Z\pi}\Z\Z^n$, and the varieties defined by them,
may be computed from the presentation matrix
$D_p\otimes_{\Z\pi}\Z\Z^n$,
where $D_p$ is the matrix given in \eqref{eq:pippres}.

\subsection{Characteristic varieties}
\label{subsec:charvar}
There is another, well-known way to associate subvarieties of the
complex algebraic torus to a space $X$ satisfying conditions
\eqref{cond:finite} and \eqref{cond:free} from \S\ref{subsec:min}.
For positive integers $k$ and $p$, set
$V^p_k(X)=\{t\in (\C^{*})^n \mid
\dim_{\C} H_p(C_*(\wX)\otimes_{\Z\pi}\C_{t}) \ge k\}$,
where $n=b_1(X)$, and $\C_{t}$ is the $\pi$-module $\C$,
given by the representation
$\pi\xrightarrow{\text{ab}}\Z^n\to\C^*$, gotten by sending
$x_i$ to $t_i$. This is an algebraic subvariety of $(\C^*)^n$,
called the {\em $(p,k)$-characteristic variety} of $X$. It is
straightforward to extend this definition to the relative
setting, as follows.

\begin{definition}
\label{def:relcv}
Let $(Y,X)$ be a CW-pair of spaces satisfying conditions
\eqref{cond:finite} and \eqref{cond:free}, and
such that the inclusion $X\hookrightarrow Y$ induces an isomorphism on
$\pi_1$.
Set $n=b_1(X)$. For $k$, $p>0$, the {\em $(p+1,k)$-characteristic
variety}
of $(Y,X)$ is the subvariety of the complex algebraic $n$-torus defined
by
\[
V^{p+1}_k(Y,X):=\{t\in (\C^{*})^n \mid
\dim_{\C} H_{p+1}(C_*(\wY,\wX)\otimes_{\Z\pi}\C_{t}) \ge k\}.
\]
\end{definition}

The characteristic variety $V_k^1(X)$ was interpreted by Hironaka
\cite{Hr}
as the variety defined by $F_{k+1}(X)$, the ideal generated
by the codimension $k$ minors of the Alexander matrix of $\pi_1(X)$:
\[
V_k^1(X)=V(F_{k+1}(X))
\]
  The next result
provides a higher-dimensional analogue of Hironaka's theorem, in a
relative setting.

\begin{prop}
\label{prop:cvpip}
Let $Y$ be a minimal $K(\pi,1)$-complex satisfying condition
\eqref{cond:surj}
from \S\ref{subsec:min}, and let $X\subset Y$ be a proper, connected
subcomplex, such that $X^{(2)}=Y^{(2)}$.
Set $p=p(X)$. Then, for all $k\ge 0$,
\[
V_k(\pi_p(X)) = V_k^{p+1}(Y,X).
\]
\end{prop}

\begin{proof}
By Remark~\ref{rem:subcomplex}, $\pi_p(X)$ is a finitely
presented $\Z\pi$-module, with presentation matrix $\Delta_p$
given in \eqref{eq:deltamin}. Hence, $V_k(\pi_p(X))=
\{t\in (\C^*)^n\mid \dim\coker \Delta_p(t)\ge k\}$,
where
$\Delta_p(t)$ is the evaluation of the matrix of Laurent polynomials
$\Delta_p\otimes_{\Z\pi}\C\Z^n$ at $x_i=t_i$.

Let $j\colon X\to Y$ be the inclusion. The
lift to universal covers, $\tilde{j}\colon \wX\to \wY$, gives rise to 
an
exact sequence of $\pi$-equivariant chain complexes,
a fragment of which is shown below:
\begin{equation}
\label{eq:ladder3}
\xymatrix @C5pc {
  H_{p+2}(Y)\otimes\Z\pi \ar@{->>}[r]^{\Pi_{p+2}\otimes \id \ }
\ar[d]^{\partial_{p+2}}
& H_{p+2}(Y,X)\otimes\Z\pi
\ar[d]^{\bar\partial_{p+2}}
\\
  H_{p+1}(Y)\otimes\Z\pi \ar@{->>}[r]^{\Pi_{p+1}\otimes \id \ }
& H_{p+1}(Y,X)\otimes\Z\pi.
}
\end{equation}
Now fix $t\in (\C^*)^n$. Tensoring \eqref{eq:ladder3}
over $\Z\pi$ with $\C$, via the representation $t\colon \pi\to \C^*$,
yields the commuting diagram
\begin{equation*}
\label{eq:ladder4}
\xymatrix  @C5pc {
H_{p+2}(Y)\otimes\C  \ar@{->>}[r]^{\Pi_{p+2}\otimes \id\ }
\ar[d]^{\partial_{p+2}(t)}
& H_{p+2}(Y,X)\otimes\C  \ar[d]^{\bar\partial_{p+2}(t)}
\\
H_{p+1}(Y)\otimes\C  \ar@{->>}[r]^{\Pi_{p+1}\otimes \id\ }
& H_{p+1}(Y,X)\otimes\C.
}
\end{equation*}
(Note also that $H_p(Y,X)=0$, by the definition of $p$.)
Chasing this diagram, we see that
$H_{p+1}(\wY,\wX;\C_t)= \coker \bar\partial_{p+2}(t)= \coker
\Delta_p(t)$.
  From Definition~\ref{def:relcv}, we have
$V^{p+1}_k(Y,X)= \{t\in (\C^{*})^n \mid
\dim H_{p+1}(\wY,\wX; \C_{t}) \ge k\}$, and we are done.
\end{proof}

\section{Higher homotopy groups of hypersolvable arrangements}
\label{sec:hypsolv}

We now apply the machinery developed in \S\ref{sec:general}
to the class of spaces we had in mind all along:
complements of complex hyperplane arrangements.
We begin with a brief review of basic notions and
relevant general results.

\subsection{Minimal cell decompositions of arrangements}
\label{subsec:minarr}
A {\em (complex) hyperplane arrangement} is a finite
set, $\A$, of codi\-mension-$1$ affine subspaces in a
finite-dimensional
complex vector space, $V$. The two main objects associated
to an arrangement $\A$ are its {\em complement},
$X(\A)=V\setminus\bigcup_{H\in \A}H$, and its
{\em intersection lattice},
$L(\A)=\left\{\bigcap_{H\in \B}\, H\mid \B\subseteq\A\right\}$.
A general reference for the subject is the book by Orlik and Terao
\cite{OT}.

Since $X(\A)$ is the complement of a complex hypersurface,
it has the homotopy type of a finite CW-complex,
and thus satisfies condition~(\ref{cond:finite}) from
\S\ref{subsec:min}.
Explicit regular CW-complexes (of dimension equal to $\dim_{\C} V$)
onto which $X(\A)$ deform-retracts were given by Salvetti \cite{Sa}
(in the complexified-real case), and by Bj\"orner and Ziegler \cite{BZ}
(in the general case). Neither of these complexes, though, is
minimal.

In a first preprint version of \cite{PS}, we proved that the
complements of arbitrary complex hyperplane arrangements satisfy
the minimality condition \eqref{eq:minimal}, up to $q=2$, by
combining results from \cite{Li1} and \cite{CSbm}. Since then,
the minimality question for arrangement complements, raised in
preprint version $1$ of \cite{PS}, has been solved in the affirmative
by Dimca~\cite{D} and Randell~\cite{R} (independently).
Using Morse theory, they proved the following.

\begin{thm}[Dimca \cite{D}, Randell \cite{R}]
\label{thm:dr}
Let $\A$ be a complex hyperplane arrangement, with
complement $X(\A)$.  Then $X(\A)$ is minimal, i.e.,
it admits a finite cell decomposition with number of $q$-cells
equal to the $q$-th Betti number, for all $q$.
\end{thm}

As noted in \cite{D}, complements of generic projective
hypersurfaces fail to possess minimal cell structures.
This indicates that minimality is a strong topological
peculiarity of complements of complex arrangements.

Assume now that $\A$ is an {\em essential $k$-generic section} of
an essential, central, aspherical arrangement $\wA$, with $k\ge 2$.
That is, $\wA$ is an essential arrangement of
hyperplanes passing through the origin of $V$, the complement
$X(\wA)$ is aspherical, and there is a complex linear subspace,
$U\subset V$, such that $U$ is $L_{k}(\wA)$-generic (in the sense of
definition \S $5(1)$ from \cite{DP}), and the restriction
$\wA^{U}$ is equal to $\A$. In this framework,
the refined Morse-theoretic analysis from \cite[Section 5]{DP}
leads to the following improved minimality result.

\begin{theorem}
\label{thm:minarr}
Let $\A$ be an essential, $k$-generic section ($k\ge 2$) of an 
essential,
aspherical arrangement $\wA$, as above. Set $p=p(X(\A))$. Then
there is a cellular classifying map between minimal 
CW-complexes with cohomology rings
generated in degree $1$, $j\colon X \to Y$, which restricts to the 
identity on
$p$-skeleta, and has the homotopy type of the inclusion,
$X(\A)\hookrightarrow X(\wA)$.
\end{theorem}
\begin{proof}
Denote by $X'(\A)$ and $X'(\wA)$ the complements of the associated
projective arrangements in the projective spaces $\P(U)$ and $\P(V)$,
respectively. The triviality of Hopf fibrations of arrangement 
complements
readily gives homotopy equivalences,
$X(\A)\simeq X'(\A) \times S^1$ and $X(\wA) \simeq X'(\wA) \times S^1$,
and implies that $p(X'(\A))=p$.

Propositions $14$ and $15$ from \cite{DP} may be used to replace,
up to homotopy, the inclusion $X'(\A)\hookrightarrow X'(\wA)$ by
a cellular map between minimal CW-complexes, $j'\colon X'\to Y'$,
which restricts to the identity on $p$-skeleta;
see the discussion preceding Theorem 16 from \cite{DP}.
Let $X=X' \times S^1$ and $Y=Y' \times S^1$.
As shown by Brieskorn \cite{Br}, all complements of
complex hyperplane arrangements satisfy
condition \eqref{cond:surj} from \S\ref{subsec:min}.
Thus, $X$ and $Y$ satisfy that condition, too.
Now set  $j:= j' \times \id\colon X\to Y$.
The claimed properties of the map $j$ follow  
from the corresponding properties of $j'$, 
and the fact that $p\ge 2$ (by Remark~\ref{rem:px}),
which guarantees that $j$ is a classifying map.
\end{proof}

\subsection{OS-algebras}
\label{subsec:osalg}
The minimality result from 
Theorem~\ref{thm:minarr} opens the way for using our approach
to generalize Hattori's results to a
wider class of arrangements.  But first,
we need to recall an important result of Orlik and
Solomon~\cite{OS}, which gives a combinatorial
description of cohomology rings of arrangement complements.

Let $\A=\{H_1,\dots, H_n\}$ be a {\em central} arrangement.
By definition, the {\em OS-algebra} of $\A$ is
\begin{equation}
\label{eq:OS}
\OS(\A) = \left.\sideset{}{^*}\bigwedge (e_1,\dots,e_n) \right\slash
\Big(
\partial e_{\B}\:\big|\, \B\subset \A\ \text{and} \
\codim \bigcap_{H\in \B} H < \abs{\B}\Big),
\end{equation}
where $\bigwedge^*(e_1,\dots,e_n)$ is the exterior
algebra over $\Z$ on generators $e_i$ in degree $1$, and
for $\B=\{H_{i_1},\dots, H_{i_r}\}$,
$e_{\B}=e_{i_1} \cdots e_{i_r}$ and
$\partial e_{\B}=\sum_{q}(-1)^{q-1}e_{i_1} \cdots
\widehat{e_{i_q}}\cdots e_{i_r}$.
There is then an isomorphism of graded algebras,
\begin{equation}
\label{eq:piqos}
\operatorname{OS}\colon H^*(X(\A))\cong\OS(\A).
\end{equation}
Under this identification, the basis $\{e_1,\dots ,e_n\}$ of
$\operatorname{A}^1(\A)$ is dual to the basis of $H_1(X(\A))$
given by the meridians of the hyperplanes, see \cite{OT}.
With respect to a fixed ordering of the hyperplanes, a canonical basis
for
$\OS(\A)$ is the {\em no broken circuits} (or, \textbf{nbc}) basis,
see \cite{OT}.

There is another, closely related, graded algebra, $\QOS(\A)$,
called the {\em quadratic Orlik-Solomon algebra},
defined as the quotient of $\bigwedge^*(e_1,\dots,e_n)$ by
relations of the form $\partial e_{\B}$, for all $\B\subset \A$
such that $\codim \bigcap_{H\in \B} H <\abs{\B}$ and $\abs{\B}=3$,
see \cite{Fa1, SY}. Clearly, the algebra $\OS(\A)$ is a quotient of
$\QOS(\A)$, and the two algebras coincide up to degree~$2$. Denote by
\begin{equation}
\label{eq:osiso}
\pi^*_{\A}\colon\QOS(\A) \to \OS(\A)
\end{equation}
the canonical projection.
Also denote by $P_{\A}(T)$ the Poincar\'e polynomial of $\OS(\A)$,
and by $\o{P}_{\A}(T)$ that of $\QOS(\A)$.
It follows at once that $\o{P}_{\A}(T) \succcurlyeq P_{\A}(T)$
(coefficientwise inequality).

\subsection{Supersolvable and hypersolvable arrangements}
\label{subsec:hyparr}

Perhaps the best understood arrangements are the supersolvable
(or, fiber-type) arrangements, introduced by Falk and Randell in
\cite{FR}.
A central arrangement $\A$ is called {\em supersolvable} if
its intersection lattice is supersolvable, in the sense of
Stanley~\cite{St}. For our purposes here, another (equivalent)
combinatorial definition will be, however, more useful;
see Definition~\ref{def:ssseries}.
The standard example is the braid arrangement in $\C^{\ell}$,
$\B_{\ell}=\{\ker (z_i-z_j)\}_{1\le i<j\le \ell}$, with
$L(\B_{\ell})=\mathcal{P}_{\ell}$, the partition lattice, and
$\pi_1(X(\B_{\ell}))=P_{\ell}$, the pure braid group on $\ell$ strings.
It follows from a theorem of Terao \cite{Te} and results in \cite{FR}
that the complement of an arbitrary
supersolvable arrangement is a $K(\pi,1)$.

The class of hypersolvable arrangements actually motivated the
framework for
  our Theorem \ref{thm:hure}. We start by reviewing the definition and
basic
properties of such arrangements.

Let $\A =\{H_{1},\dots ,H_{n}\}$ be a central arrangement in the
complex
vector space $V$. Denote also by $\A=\{\a_1,\dots ,\a_n\}\subset
\P(V^*)$
its set of defining equations, viewed as points in the dual projective
space.
Let $\B\subset \A$ be a proper, non-empty sub-arrangement, and set
$\o{\B}:=\A\setminus \B$. We say that $(\A,\B)$ is a {\em solvable
extension} if the following conditions are satisfied (see \cite{JP1}):

\renewcommand{\thesubfigure}{(\Roman{subfigure})}

\begin{figure}[ht]
\setlength{\unitlength}{0.5cm}
\subfigure[]{
\label{fig:ax1}%
\begin{minipage}[t]{0.23\textwidth}
\begin{picture}(4,4)(-1.5,0)
\put(0.5,0.5){\circle*{0.4}}
\put(2,2){\circle*{0.4}}
\put(3.5,3.5){\circle{0.4}}
\put(0.6,0.6){\line(1,1){1.25}}
\put(2.1,2.1){\line(1,1){1.25}}
\put(1.2,0.5){\makebox(0,0){$\alpha$}}
\put(2.7,2){\makebox(0,0){$\beta$}}
\put(4.2,3.5){\makebox(0.1,0){$\not\!\exists\, a$}}
\end{picture}
\end{minipage}
}
\subfigure[]{
\label{fig:ax2}%
\begin{minipage}[t]{0.28\textwidth}
\begin{picture}(4,4)(-1.5,0)
\put(0.5,0.5){\circle{0.4}}
\put(2,2){\circle{0.4}}
\put(3.5,3.5){\circle*{0.4}}
\put(0.65,0.65){\line(1,1){1.2}}
\put(2.15,2.15){\line(1,1){1.25}}
\put(1.2,0.5){\makebox(0,0){$a$}}
\put(2.7,2){\makebox(0,0){$b$}}
\put(4.2,3.5){\makebox(2.3,0){$\exists ! f(a,b)$}}
\end{picture}
\end{minipage}
}
\subfigure[]{
\label{fig:ax3}%
\begin{minipage}[t]{0.28\textwidth}
\begin{picture}(6,4)(-1,0)
\put(0.5,0.5){\circle{0.4}}
\put(2,2){\circle{0.4}}
\put(3.5,3.5){\circle*{0.4}}
\put(3.5,1.5){\circle{0.4}}
\put(5,2){\circle*{0.4}}
\put(6.5,0.5){\circle*{0.4}}
\put(0.65,0.65){\line(1,1){1.2}}
\put(2.15,2.15){\line(1,1){1.25}}
\put(3.5,3.5){\line(1,-1){3}}
\put(2.2,1.92){\line(3,-1){1.09}}
\put(0.71,0.55){\line(3,1){2.59}}
\put(4.8,1.92){\line(-3,-1){1.09}}
\put(6.29,0.55){\line(-3,1){2.59}}
\put(1.15,0.25){\makebox(0,0){$a$}}
\put(2.7,2.2){\makebox(0,0){$b$}}
\put(3.5,1){\makebox(0,0){$c$}}
\put(4.2,3.6){\makebox(1.4,0){$f(a,b)$}}
\put(5.8,2.2){\makebox(1.4,0){$f(a,c)$}}
\put(7.2,0.3){\makebox(1.4,0){$f(b,c)$}}
\end{picture}
\end{minipage}
}
\caption{\textsf{Axioms for solvable extensions}}
\label{fig:axioms}
\end{figure}

\begin{Romenum}
\item \label{I}
No point $a\in \o{\B}$ sits on a projective line determined
by $\a, \beta\in\B$.
\item \label{II}
For every $a, b\in\o{\B}$, $a\ne b$, there exists a point $\a\in \B$ on
the line passing through $a$ and $b$. (In the presence of
condition~\eqref{I}, this point is uniquely determined, and
will be denoted by $f(a,b)$.)
\item \label{III}
For every distinct points $a, b, c\in\o{\B}$, the three points
$f(a,b)$, $f(a,c)$, and $f(b,c)$ are either equal or collinear.
\end{Romenum}

Note that only two possibilities may occur:
either $\rank(\A)=\rank(\B)+1$ ({\em fibered case}), or
$\rank(\A)=\rank(\B)$ ({\em singular case});
see \cite[Lemma 1.3(i)]{JP1}.

\begin{definition}[\cite{JP1}]
\label{def:hyper}
The arrangement $\A$ is called {\em hypersolvable} if it has a
{\em hypersolvable composition series}, i.e., an ascending chain
of sub-arrangements, $\A_1\subset \cdots \subset \A_i\subset
\A_{i+1}\subset \cdots \subset\A_{\ell}=\A$, where $\rank \A_1=1$, and
each extension $(\A_{i+1},\A_{i})$ is solvable.
\end{definition}

The length of a composition series depends only on $\A$; it will
be denoted by $\ell(\A)$. Note that the property of being hypersolvable
is purely combinatorial. In fact, given an arrangement $\A$, one can
decide whether it is hypersolvable or not, only from the
elements of rank one and two of $L(\A)$, since the definitions
only involve the collinearity relations in $\A$.
The class of hypersolvable arrangements
includes supersolvable arrangements, cones of generic arrangements
(for which $\ell(\A)= \abs{\A}$),
and many others, see \cite{JP1}, and the examples
in \S\S\ref{subsec:randell},
\ref{subsec:pi2structure}, and \ref{subsec:notriangle}.

It is appropriate to mention here that one may assume from now on,
whenever necessary, that a given hypersolvable arrangement $\A$
is also essential. Indeed, one knows how to associate to
an arbitrary central arrangement $\A$ an essential arrangement,
$\A_{\rm ess}$, without changing the homotopy types of the complements,
$X(\A)$ and $X'(\A)$, and the intersection lattice $L(\A)$; see
\cite[p.~197]{OT}.

The connection between hypersolvable and supersolvable (or fiber-type)
arrangements comes from the following fact, which is implicit in
\cite[Lemma 4.5]{JP1}, and is explicitly proved in
\cite[Prop.~1.3(i)]{JP2}.
If the solvable extension $(\A, \B)$ is fibered,
then there is a Serre fibration $X(\A) \to X(\B)$,
with homotopy fiber $\C \setminus \{\text{$m$ points} \}$,
where $m=|\overline{\B}|$. It follows
from \cite[Prop.~1.3]{JP2} that the (topological)
definition of fiber-type arrangements may be rephrased
in hypersolvable terms, as follows.

\begin{definition}[\cite{FR}]
\label{def:ssseries}
The arrangement $\A$ is {\em supersolvable} (or, {\em fiber-type}) if
it has a {\em supersolvable composition series}, that is, a 
hypersolvable
composition series as in Definition~\ref{def:hyper}, for which all
extensions are fibered.
\end{definition}

We thus see that all fiber-type arrangements are hypersolvable. On the
other hand, one knows from \cite[Thm.~D]{JP1} that a hypersolvable
arrangement $\A$ is aspherical if and only if $\A$ is fiber-type, which
happens precisely when $\ell(\A)=\rank(\A)$.

\subsection{Supersolvable deformations}
\label{subsec:ssdef}
Our basic tool for the topological study of hypersolvable arrangements
is the following theorem, which mostly puts together and organizes
a number of known results.

\begin{theorem}
\label{thm:hypsolv}
Let $\A$ be an essential, hypersolvable arrangement, with composition
series $\A_1\subset \cdots \subset\A_{\ell}=\A$, and
{\em exponents} $d_i:=\abs{\A_i\setminus \A_{i-1}}$.
Set $\pi= \pi_1(X(\A))$ and $p=p(X(\A))$.  Then
$\A$ is a $2$-generic section of an essential,
supersolvable arrangement $\wA$, called the
{\em supersolvable deformation} of $\A$, such that:
\begin{enumerate}

\item \label{part1}
$\wA$ has a supersolvable composition series,
$\wA_1\subset \cdots \subset\wA_{\ell}=\wA$, with $|\wA_i|=\abs{\A_i}$,
for $1\le i\le \ell$.
\item \label{part2}
$X(\wA)$ sits atop a tower,
$X(\wA_{\ell})\xrightarrow{p_{\ell}}
X(\wA_{\ell-1})\to\cdots\to X(\wA_2)
\xrightarrow{p_2} X(\wA_1)=\C^*$,
of Serre fibrations, $p_i\colon X(\wA_i)\to X(\wA_{i-1})$,
with fiber $\C\setminus \{\text{$d_i$ points}\}$, and monodromy
$\rho_{i-1}\colon \pi_1(X(\wA_{i-1}))\to P_{d_i}\subset \Aut(F_{d_i})$.
\item \label{part3}
$X(\wA)$ is a $K(\pi,1)$ space.
The fundamental group
admits an iterated semidirect product decomposition,
$\pi=F_{d_{\ell}}\rtimes_{\rho_{\ell-1}}\cdots
\rtimes_{\rho_1}F_{d_1}$.
\item \label{part5}
There is a canonical isomorphism,
$\o{\operatorname{OS}}\colon H^*(X(\wA))\cong\QOS(\A)$.
Under this isomorphism, and the isomorphism
$\operatorname{OS}\colon H^*(X(\A))\cong\OS(\A)$,
the map $J^*\colon H^*(X(\wA))\to H^*(X(\A))$,
induced by the canonical inclusion, 
$J \colon X(\A)\hookrightarrow X(\wA)$,
corresponds to the canonical projection
$\pi^*_{\A}\colon \QOS(\A)\to \OS(\A)$.
\item \label{part6}
$\o{P}_{\A}(T)=P_{\wA}(T)=\prod_{i=1}^{\ell}(1+d_iT)$.
\item \label{part0}
There exist minimal CW-complexes $X$ and $Y$,
and homotopy equivalences $\phi\colon X\to X(\A)$ and
$\psi\colon X(\wA)\to Y$, such that the composite $\psi J \phi$
is homotopic to a cellular classifying map,
\[
j\colon X \rightarrow Y ,
\]
which restricts to the identity on $p$-skeleta.
\item \label{part5half}
The dual of the split exact sequence
$$0\to H_*(X)\xrightarrow{j_*}H_*(Y)\xrightarrow{\Pi_*}H_*(Y,X)\to 0$$
from \eqref{eq:jpi} may be identified with
\begin{equation}
\label{eq:dualpr}
0\to \ker(\pi^*_{\A}) \xrightarrow{\iota^*} \QOS(\A)
\xrightarrow{\pi^*_{\A}}
\OS(\A) \to 0\, ,
\end{equation}
where $\iota^*$ denotes the inclusion of the kernel.
\end{enumerate}
\end{theorem}

\begin{proof}
\eqref{part1}
The supersolvable arrangement $\wA$ is obtained from $\A$ by
the deformation method introduced in \cite{JP1}, and refined in 
\cite{JP2}.
Part~\eqref{part1} follows from this deformation method, which
proceeds inductively, using the given composition series of $\A$.

\eqref{part2} Up to homotopy, we may view each $\wA_i$ as an
arrangement in $\C^i$, and replace each map $p_i$ by a bundle map,
$q_i$, with the specified fiber (more precisely, by a linear fibration,
admitting a section, see \cite{FR}). Moreover, the defining polynomials
for $\wA_i$ may be written inductively as $f_1(z_1)=z_1$, and
$f_i(z_1,\dots,z_i)=f_{i-1}(z_1,\dots,z_{i-1})\cdot
\prod_{k=1}^{d_i}(z_i-g_{i,k}(z_1,\dots,z_{i-1}))$.
Clearly, $f_i/f_{i-1}$ is a completely solvable Weierstrass polynomial
over $X(\wA_{i-1})$. Thus, by \cite[Thm.~2.3]{CSbm}, the monodromy of
the bundle map $q_i$ factors through the pure braid group $P_{d_{i}}$,
acting on the free group $F_{d_{i}}$ via the Artin representation.

\eqref{part3} The first assertion follows from \cite{JP2}. 
The specified structure of
$\pi$ is provided by Part~\eqref{part2}.

\eqref{part5} Since $\wA$ is supersolvable,
$\OS(\wA)\cong\QOS(\wA)$, see Falk~\cite{Fa1} and Shelton and
Yuzvinsky~\cite{SY}.
Moreover, Theorem~2.4 from \cite{JP2} insures that $\A$ and $\wA$ have
the same collinearity relations, which implies that
$\QOS(\A)\cong\QOS(\wA)$. The canonical
isomorphism in \eqref{part5} is then given by:
\begin{equation}
\label{eq:homhyp}
H^*(X(\wA))\cong\OS(\wA)\cong\QOS(\wA)\cong\QOS(\A).
\end{equation}
The identification of $J^*$ with $\pi^*_{\A}$ follows from the fact
that the basis $\{e_1,\dots,e_n\}$ of $\operatorname{A}^1$ is dual
to the basis of $H_1$ given by the meridians, and
$J_1\colon H_1(X(\A))\to H_1(X(\wA))$ preserves those meridians.

\eqref{part6} The equality between the Poincar\'{e} polynomials of
$\QOS(\A)$ and $\OS(\wA)$ follows from \eqref{part5}. The second
equality follows from \cite{FR}, via Part~\eqref{part1}.

\eqref{part0}
This follows from Theorem~\ref{thm:minarr}.

\eqref{part5half}
This follows from Parts~\eqref{part5} and \eqref{part0}.
\end{proof}

We record as a corollary the most important (for our purposes)
consequence of the above theorem.

\begin{corollary}
\label{cor:hsclass}
All
hypersolvable complements, and their $K(\pi, 1)$ spaces, are minimal,
with cohomology algebra generated in degree one.
\end{corollary}

The following corollary shows that the order of $\pi_1$-connectivity of
the
complement of a hypersolvable arrangement is combinatorially determined
(though $\pi_1$ itself is not {\it a priori} combinatorial).

\begin{corollary}
\label{cor:pcomb}
Let $\A$ be a hypersolvable arrangement. Set $X=X(\A)$, and
$\pi=\pi_1(X)$.
Let $p=p(X)$ be the order of $\pi_1$-connectivity of $X$. Then:
\begin{enumerate}
\item \label{i1}
$p(X)=\sup\, \{k \mid P_{\A}(T)\equiv\o{P}_{\A}(T)\bmod (T^{k+1})\}$.
\item \label{i2}
$p(X)\ge 2$.
\item \label{i3}
$p(X)=\infty \Longleftrightarrow P_{\A}(T)=\o{P}_{\A}(T)
\Longleftrightarrow \text{$\A$ is supersolvable}$.
\item \label{i4}
If $p(X)<\infty$, then $\o{P}_{\A}(T) - P_{\A}(T) \equiv
c_{p+1}T^{p+1} \bmod (T^{p+2})$,
where $c_{p+1}$ is a {\em positive} integer.
\end{enumerate}
\end{corollary}

\begin{proof}
\eqref{i1} Follows from Theorem~\ref{thm:hypsolv},
Parts~\eqref{part3} and \eqref{part5}.

\eqref{i2} Follows from Remark~\ref{rem:px} and
Theorem~\ref{thm:hypsolv}, Part~\eqref{part0}.

\eqref{i3} Follows from \eqref{i1} and \cite[Prop.~3.4]{JP2}.

\eqref{i4} Follows from \eqref{i1} and the fact that
$\o{P}_{\A}(T) \succcurlyeq P_{\A}(T)$.
\end{proof}

\subsection{A presentation for $\pi_p(X(\A))$}
\label{subsec:firsths}
We come now to the main result in this section.
Together with Corollary~\ref{cor:pcomb}\eqref{i3}--\eqref{i4},
this result provides a complete proof of
Theorem~\ref{thm:intro1} from the Introduction.

Let $\A$ be a hypersolvable arrangement.
Set $\pi= \pi_1(X(\A))$ and $p=p(X(\A))$. Denote by $\wA$ the
supersolvable deformation of the associated essential
hypersolvable arrangement, $\A_{\rm ess}$.
Use Theorem~\ref{thm:hypsolv}\eqref{part0} to replace,
up to homotopy, the inclusion
$J \colon X(\A_{\rm ess}) \hookrightarrow X(\wA)$ by a map,
$j \colon X \to Y$, satisfying the hypotheses of 
Theorem~\ref{thm:hure}.

\begin{theorem}
\label{thm:hspipres}
Let $\A$ be a hypersolvable arrangement, with
fundamental group $\pi=\pi_1(X(\A))$ and order of
$\pi_1$-connectivity $p=p(X(\A))$.
Then:
\begin{enumerate}
\item \label{p1}
$X(\A)$ is aspherical $\Longleftrightarrow p= \infty$.
\item \label{p2}
If $p< \infty$, then the first non-vanishing higher homotopy group of
$X(\A)$ is $\pi_p(X(\A))$; as a $\Z\pi$-module, $\pi_p(X(\A))$ is
isomorphic to $\coker (D_p)$, where the presentation matrix $D_p$
is given by \eqref{eq:pippres}.
\item \label{p3}
If $p< \infty$, then the group of $\pi$-coinvariants
of $\pi_p(X(\A))$ is free abelian, of rank
\[
c_{p+1}=\text{coefficient of $T^{p+1}$ in $\o{P}_{\A}(T) - P_{\A}(T)$}.
\]
In particular, both $p$ and the group
$(\pi_p(X(\A)))_{\pi}$ are combinatorially determined.
\end{enumerate}
\end{theorem}

\begin{proof}
Parts \eqref{p1} and \eqref{p2} follow from Theorem~\ref{thm:hure} and
Corollary~\ref{cor:coinv}. Part \eqref{p3} follows from
Corollary~\ref{cor:coinv}, Theorem~\ref{thm:hypsolv}\eqref{part5half}
and Corollary~\ref{cor:pcomb}\eqref{i1}.
\end{proof}

\begin{remark}
\label{rem:matrices}
$(i)$ Let $(C_*(\wY), \partial_*)$ be the $\pi$-equivariant
chain complex of the universal cover, associated to the
minimal Morse-theoretic cell structure, $Y$, of $X(\wA)$, as in
Theorem~\ref{thm:hspipres}. Since $X(\wA)$ is aspherical,
$(C_*(\wY), \partial_*)\xrightarrow{\epsilon}\Z$ is a finite, free,
$\epsilon$-minimal, $\Z\pi$-resolution of $\Z$. In particular, it is
chain homotopy equivalent to the $\epsilon$-minimal, finite, free,
$\Z\pi$-resolution of $\Z$ constructed by Fox calculus in \cite{CScc},
starting from the iterated semidirect product structure of $\pi$,
with trivial monodromy actions on homology, described in
Theorem~\ref{thm:hypsolv}, Parts~\eqref{part2} and \eqref{part3}.

$(ii)$ More generally, let $j \colon X\to Y$ be a classifying map
which satisfies the requirements from Theorem~\ref{thm:hure}. We
can show that the second nilpotent quotient of the first higher
non-trivial homotopy group of $X$, $\pi_p(X)/I^2\, \pi_p(X)$,
is determined, as a filtered module over $\Z\pi_1(X)/I^2$, by
the map induced by $j$ between cohomology rings. For arbitrary
hypersolvable complements, the above nilpotent quotient turns out,
in this way, to be combinatorially determined. Details will appear
elsewhere.
\end{remark}

\subsection{Comparison with some results of Randell}
\label{subsec:randell}
A formula for the coinvariants of the first non-vanishing higher
homotopy group, similar to our \ref{thm:hspipres}\eqref{p3},
was obtained by Randell, using different methods, in
\cite[Thm.~2 and Prop.~9]{Ra2}, for the class of iterated generic
hyperplane sections (of rank $\ge 3$) of essential,
aspherical arrangements. For an arrangement $\A$ in this class,
$p(X(\A))=\rank(\A)-1$, by results from \cite{Ra2}.
Randell's class of arrangements, and the class of hypersolvable
arrangements have a similar behavior, from the point of view of
the coinvariants of the first higher non-vanishing homotopy group.
Nevertheless, the two classes are distinct, as the next two examples
show:

\begin{example}
\label{ex:ran1}
  For $\ell \ge 5$, let
$\A_{\ell} := \B_{\ell} \cup \{H\}$,
where $\B_{\ell}=\{z_i - z_j = 0\}_{1 \leq i < j \leq \ell}$ and
$H= \{ z_1+z_2+z_3-3z_{\ell}=0 \}$.
Each arrangement $\A_{\ell}$ is hypersolvable, of rank $\ell-1$ and
length
$\ell$. We claim that
$p(X(\A_{\ell})) = 2$. It follows that these arrangements cannot be
iterated generic sections of essential, aspherical arrangements,
since
this would imply that $p(X(\A_{\ell})) ={\ell}-2$.

The claim may be verified by showing that
$\rank A^3(\A_{\ell})< \rank \overline{A}^{3}(\A_{\ell})$;
see Corollary~\ref{cor:pcomb}, Parts~\eqref{i1} and \eqref{i2}.
Let $\mathcal{C}= \{H_1,H_2,H_3,H \}$, with $H_i= \{z_i-z_{\ell}=0 \}$,
and let $\{e_1,e_2,e_3,e \}$ be the corresponding OS-generators.
It is easy to check that
$\rank A^3(\A_{\ell}) \leq \rank
\left(\overline{A}(\B_{\ell}) \otimes
\bigwedge^*(e)/(\partial e_{\mathcal{C}})\right)^3$, and
$\rank \overline{A}^3(\A_{\ell}) =\rank \left(\overline{A}(\B_{\ell})
\otimes
\bigwedge^*(e)\right)^3$,
directly from the definitions (see \S\ref{subsec:osalg}).
Notice that $\{ H_1,H_2,H_3 \}\subset \B_{\ell}$ is a boolean
subarrangement,
hence $e_1e_2e_3$ is a non-zero element of $\overline{A}(\B_{\ell})$
(use \cite[Prop.~3.66]{OT}). We infer that $\partial e_{\mathcal{C}}$
is a
non-zero element of $\overline{A}(\B_{\ell})\otimes \bigwedge^*(e)$,
whence the desired inequality.
\end{example}

\begin{example}
\label{ex:ran2}
Let $\A$ be an iterated generic section of an essential,
aspherical arrangement $\B$ which is not hypersolvable.
For example, take $\B$ to be the reflection arrangement
of type $\operatorname{D}_n$, with $n \geq 4$, see \cite{JP1}.
If $\rank (\A) \geq 3$, then necessarily $\A$ and $\B$
have the same collinearity relations, and therefore $\A$
cannot be hypersolvable.
\end{example}

There is, however, a certain overlap between the two classes.
For instance, iterated generic sections (of rank $\geq 3$) of
fiber-type arrangements are obviously hypersolvable.

\begin{example}
\label{ex:overlap}
Let $\A$ be an essential, proper,
$k$-generic section of an essential arrangement $\wA$, with
$k= \rank (\A)-1$. It follows that $\A$ must be an iterated
generic hyperplane section of $\wA$; see the proof of Proposition $14$
from \cite{DP}. As a particular case, consider a hypersolvable 
arrangement
$\A$ in $\C^3$, such that $p(X(\A))< \infty$ (i.e., such that
$\A$ is not supersolvable, see Corollary~\ref{cor:pcomb}\eqref{i3}).
We must have $\rank (\A)=3$, since all arrangements are fiber-type
in rank $\le 2$. Consequently, $\A$ is a $2$-generic, proper
section of its essential supersolvable deformation, $\wA$, \cite{JP2},
hence also an iterated generic section of the aspherical arrangement 
$\wA$.
\end{example}

\section{On the structure of $\pi_2$ as a $\Z\pi_1$-module}
\label{sec:pi2pres}

We now analyze in more detail
the structure of $\pi_2(X')$, viewed as a module over $\Z\pi_1(X')$,
in the case when $X'=X(\A')$ is the complement of an affine line
arrangement
whose cone $\A=\bf{c}\A'$ is hypersolvable.

\subsection{$K(\pi,1)$ tests}
\label{subsec:kpi1test}

Let $\A'$ be an arrangement of affine lines in $\C^2$.
The complement $X'=X(\A')$ has the homotopy-type of a
$2$-complex, hence the only obstruction to $X'$ being
aspherical is the second homotopy group, $\pi_2(X')$.
In \cite{Fa2}, Falk gave several conditions (some
sufficient, some necessary), for the vanishing of
$\pi_2(X')$, providing a (partial) $K(\pi,1)$-test
for complexified line arrangements. This test is
geometric in nature, involving Gersten-Stallings
weight systems.

Another partial $K(\pi, 1)$-test,
valid this time in all dimensions, but only for
hypersolvable arrangements, was given in \cite[Thm.~D]{JP1}.
This test is purely combinatorial. Assuming the cone
$\A=\bf{c}\A'$ is hypersolvable, it says that $X'$ is
aspherical if and only if $\ell(\A)= \rank (\A)$.

None of these asphericity tests, though, gives a precise description
of $\pi_2(X')$, viewed as a $\pi_1(X')$-module. Our machinery
affords such a description, at least in the special case when $\A$
is hypersolvable.

\subsection{Affine arrangements with hypersolvable cones}
\label{subsec:affarr}
We first describe the structure of the fundamental
group of the complement of a deconed fiber-type arrangement,
of arbitrary rank.
\begin{lemma}
\label{lem:decone}
Let $\A$ be a supersolvable arrangement, with
composition series $\A_1\subset \cdots \subset \A_{\ell}$, and let
$F_{d_{\ell}}\rtimes_{\rho_{\ell-1}}
F_{d_{\ell-1}}\rtimes\cdots\rtimes_{\rho_2} F_{d_2}
\rtimes_{\rho_1} F_1$ be the corresponding iterated semidirect
product decomposition of
$\pi_1(X(\A))$. If ${\bf{d}}\A$ is a decone of $\A$, then
$\pi_1(X({\bf{d}}\A))=F_{d_{\ell}}\rtimes_{\rho_{\ell-1}}\cdots
\rtimes_{\rho_2} F_{d_2}$.
\end{lemma}

\begin{proof}
Recall from the proof of Theorem~\ref{thm:hypsolv}\eqref{part2}
that $\A$ has defining polynomial of the form
$f_{\A}=f_1f_2\cdots f_{\ell}$, where $f_1(z)=z_1$, and
$f_i/f_{i-1}$ is a completely solvable Weierstrass polynomial
over $X(\A_{i-1})$. The decone ${\bf{d}}\A$,
obtained by setting $z_1=1$, has defining polynomial
$f_{{\bf{d}}\A}(z_2,\dots,z_{\ell})=f_2(1,z_2)\cdots
f_{\ell}(1,z_2,\dots,z_{\ell})$.
The result follows at once.
\end{proof}

The above result gives the structure of fundamental groups
of complements of deconed hypersolvable arrangements
in all dimensions, via the deformation method from
Theorem~\ref{thm:hypsolv}, Parts~\eqref{part1}--\eqref{part3}
and \eqref{part0}.
In dimension $3$, this can be much improved, to a precise
description of the homotopy type and of the second homotopy group,
as follows.

\begin{theorem}
\label{thm:pi2hyper}
Let $\A'$ be an affine line arrangement in $\C^{2}$,
such that $\A=\bf{c}\A'$ is hypersolvable.
Set $\ell=\ell(\A)$, $X'=X(\A')$, $\pi'=\pi_1(X')$, and $p= p(X')$.
Denote by $\A_{\rm ess}$ the associated essential hypersolvable
arrangement, with supersolvable deformation $\wA$, as in
Theorem~\ref{thm:hypsolv}. Then $X'$ has the homotopy type
of the $p$-skeleton, $Y'^{(p)}$, of a minimal cell structure, $Y'$, for
$X(\bf{d} \wA)$. In particular:
\begin{enumerate}
\item \label{r1}
$X'$ is aspherical $\Longleftrightarrow \ell \leq 3 \Longleftrightarrow
p \neq 2.$
\item \label{r2}
If $\ell>3$, then $\pi_2(X')$ is non-trivial, and admits the following
finite, free, $\epsilon$-minimal $\Z\pi'$-resolution:
\begin{equation}
\label{eq:pi2res}
0\to C_{\ell-1}\xrightarrow{\partial_{\ell-1}}C_{\ell-2}\to\cdots\to
C_4\xrightarrow{\partial_{4}}C_{3}\to \pi_2(X')\to 0,
\end{equation}
where $(C_{*},\partial_{*})$ is the $\pi'$-equivariant chain complex of
$\widetilde{Y'}$.
\end{enumerate}
\end{theorem}

\begin{proof}
If $p= \infty$, then $\wA= \A_{\rm ess}$, by 
Corollary~\ref{cor:pcomb}\eqref{i3}.
Up to homotopy, $X'\simeq X(\bf{d} \wA)$, and we may use
the Morse-theoretic minimal structure from \cite[Corollary 6]{DP}.
If $p< \infty$, we know from Example~\ref{ex:overlap} that
$\A= \A_{\rm ess}$ is an iterated generic hyperplane section of
$\wA$, with $p=2$; see also the first paragraph of 
\S\ref{subsec:randell}.
The method of proof of Corollary $6$ from \cite{DP} provides the
desired minimal complex $Y'$.

\eqref{r1}
The space $X'$ is aspherical if and only if $\ell \leq 3$, by
\cite[Thm.~D]{JP1} (since $\A$ is fiber-type, if $\ell \leq 3$).
If $p \neq 2$, then necessarily $p>2$ (by Remark~\ref{rem:px}).
It follows from Theorem~\ref{thm:hspipres}\eqref{p2} that
$\pi_2(X')=\pi_2(X(\A))=0$, and so $X'$ must be
aspherical. Conversely, if $p=2$, then $\pi_2(X')$ must be non-zero  
(use
Corollary~\ref{cor:coinv}).

\eqref{r2} If $\ell>3$, we know from Part~\eqref{r1} that
$X'\simeq Y'^{(2)}$, where $\dim Y'= \ell -1$. Everything
then follows from Remarks~\ref{rem:hatres} and \ref{rem:subcomplex}.
\end{proof}

\begin{remark}
\label{rem:manyell}
The resolution \eqref{eq:pi2res} may have arbitrary length.
Indeed, for each $\ell\ge 1$, there exists a hypersolvable
arrangement $\A$ in $\C^3$ with $\ell(\A)=\ell$, see \cite[\S1]{JP1}.
Moreover, any sequence of exponents,
$\{ 1=d_1,d_2, \dots, d_{\ell} \}$, may be obtained in this way.
\end{remark}

\subsection{Structure of $\pi_2$ of a hypersolvable
line arrangement complement}
\label{subsec:pi2structure}

The group of $\pi'$-coinvariants of $\pi_2(X')$ is very simple
to describe: By Theorem \ref{thm:pi2hyper}\eqref{r2},
it is free abelian, of rank $b_3(X(\mathbf{d} \wA))=b_3(\pi')$.
On the other hand, the following result shows that $\pi_2 (X')$,
when non-trivial,
has a fairly complicated structure as a $\Z\pi'$-module.

\begin{theorem}
\label{thm:free}
Let $\A'$ be an affine line arrangement in $\C^{2}$
such that $\A=\bf{c}\A'$ is hypersolvable. Let $\ell$ be
the length of $\A$, and $\{1=d_1,d_2,\dots ,d_{\ell}\}$ the exponents.
Set $X'=X(\A')$ and $\pi'=\pi_1(X')$. Assume $\ell>3$
(so that $\pi_2(X')\ne 0$). Then:
\begin{enumerate}
\item \label{s1}
$\pi_2(X')$ is a projective $\Z\pi'$-module if and only if $\ell=4$.
In that case, $\pi_2(X')$ is free, with
rank equal to $b_3(\pi')=d_{2}d_{3}d_{4}$.
\item \label{s2}
$\pi_2 (X')$ is neither finitely generated as an abelian group,
nor nilpotent as a $\Z\pi'$-module.
\end{enumerate}
\end{theorem}

\begin{proof} \eqref{s1} From resolution~\eqref{eq:pi2res}, we see that
$\pi_2(X')$ is isomorphic to $\coker(\partial_4) = \im (\partial_3)
\subset C_2$.
If $\ell=4$, then $\pi_2(X')=C_3$ is a free $\Z\pi'$-module,
with rank $b_3(\pi')$ given by Theorem~\ref{thm:hypsolv}\eqref{part6}.
If $\ell>4$, then $\pi_2(X')$ is not projective, by the minimality of
\eqref{eq:pi2res}.

\eqref{s2} Note first that the $I$-adic filtration of the
group algebra $\Q\pi'$ is Hausdorff, in the sense that $\bigcap_{k\ge
0}
I^k =0$, where $I=\ker(\epsilon\colon\Q\pi'\to\Q)$ is the
augmentation ideal.
This follows from the fact that $\pi'$ is an iterated semidirect
product
of free groups, where all homology monodromy actions are trivial
(cf.~Lemma~\ref{lem:decone} and
Theorem~\ref{thm:hypsolv}, Parts~\eqref{part2} and \eqref{part3});
therefore, $\pi'$ is residually torsion-free
nilpotent (see \cite{FR2}), and so the $I$-adic filtration of $\Q\pi'$
must be Hausdorff (see \cite{Ch}). It follows that the $I$-adic
filtration of the
free $\Q\pi'$-module $C_2 \otimes \Q$ is also Hausdorff.

Assume now that either $\pi_2(X')$ is finitely generated as
an abelian group, or nilpotent as a $\Z\pi'$-module.
It follows that $I^k \cdot \pi_2(X') \otimes \Q = 0$, for
some $k\ge 0$, and thus $\pi_2 (X')\otimes \Q $ must be a nilpotent,
non-trivial $\Q\pi'$-module. This implies that
$(g_1-1) \cdots (g_k-1) \cdot b = 0$,
for some $g_1, \dots, g_k \in \pi' \setminus \{1\}$, and
$b \in \Q\pi' \setminus \{0\}$. On the other hand,
$\Q\pi'$ has no zero-divisors, since $\pi'$ is residually torsion-free
nilpotent (see \cite{Pa}). This gives the desired contradiction,
proving \eqref{s2}.
\end{proof}

\renewcommand{\thesubfigure}{(\alph{subfigure})}

\begin{figure}[ht]
\setlength{\unitlength}{0.6cm}
\subfigure[]{
\label{fig:ex1}%
\begin{minipage}[t]{0.25\textwidth}
\begin{picture}(3,2.8)(0.5,0.5)
\multiput(1.5,1)(0,2){2}{\line(1,0){4}}
\multiput(2.5,0)(2,0){2}{\line(0,1){4}}
\put(1.5,0){\line(1,1){4}}
\end{picture}
\end{minipage}
}
\subfigure[]{
\label{fig:ex2}%
\begin{minipage}[t]{0.25\textwidth}
\begin{picture}(3,2.8)(0.5,0.5)
\multiput(1.5,1)(0,2){2}{\line(1,0){4}}
\multiput(2.5,0)(2,0){2}{\line(0,1){4}}
\put(1.5,0.5){\line(1,1){3.5}}
\put(2,0){\line(1,1){3.5}}
\end{picture}
\end{minipage}
}
\subfigure[]{
\label{fig:ex3}%
\begin{minipage}[t]{0.25\textwidth}
\begin{picture}(3,2.8)(0.5,0.5)
\put(1,1){\line(1,0){5.5}}
\multiput(3,0)(1.5,0){2}{\line(0,1){4}}
\put(1,0){\line(1,1){4}}
\put(2.5,4){\line(1,-1){4}}
\end{picture}
\end{minipage}
}
\caption{\textsf{Some line arrangements whose
cones are hypersolvable}}
\label{fig:linearr}
\end{figure}

\begin{example}
\label{ex:1}
Let $\A'$ be the affine line arrangement from Figure~\ref{fig:ex1},
with defining polynomial $f_{\A'}=z_1z_2(z_1-1)(z_2-1)(z_2-z_1)$.
Then $\A=\mathbf{c}\A'$ is an essential $3$-slice of the braid
arrangement $\B_4$.
Hence, $\A$ is supersolvable, with length $\ell=3$, and exponents
$\{1,2,3\}$.
We then have $\pi'=F_3\rtimes F_2$, and $X'=K(\pi',1)$.
\end{example}

\begin{example}
\label{ex:2}
Let $\A'$ be the arrangement from Figure~\ref{fig:ex2}, with defining
polynomial
$f_{\A'}=(z_1-1)(z_1+1)(2z_1-2z_2-1)(2z_1-2z_2+1)(3z_1-6z_2-1)(3z_1-6z_2+1)
$.
Then $\A=\mathbf{c}\A'$ is the arrangement from Fan \cite[\S3.I]{Fn}.
It is readily seen that $\A$ is hypersolvable, with $\ell=4$,
and $\exp(\A)=\{1,2,2,2\}$. We then have $\pi'=F_2\times F_2 \times
F_2$,
and $\pi_2(X')=(\Z\pi')^8$. Notice that $V_1(\pi_2(X'))=(\C^*)^6$.
\end{example}

\begin{example}
\label{ex:3}
Let $\A'$ be the arrangement from Figure~\ref{fig:ex3}, with defining
polynomial $f_{\A'}=z_1z_2(z_1-1)(z_2-z_1-1)(z_2+z_1-2)$. Then
$\A=\mathbf{c}\A'$ is hypersolvable, with $\ell=5$, and
$\exp(\A)=\{1,1,1,1,2\}$. The algorithm from \cite{CSbm} yields
a ``braid monodromy" presentation for $\pi'$, with
generators  $x_1,\dots , x_5$, and commutation relations
$[x_i,x_j]=1$,  for $1\le i\le 3$ and $i<j$.
By \cite{Li1}, $X'$ has the homotopy type of the $2$-complex
associated to this presentation. In turn, this $2$-complex
is homotopy equivalent to the $2$-skeleton of $Y'= T^{3}\times (S^1 
\bigvee S^1)$,
the product of the $3$-torus with a wedge of two circles,
endowed with the standard minimal cell structure. From
Remarks~\ref{rem:hatres} and \ref{rem:subcomplex}, we
infer that $\pi_2(X')= \coker (\partial_4)$, where
$(C_*, \partial_*)$ is the $\pi'$-equivariant chain complex
of $\widetilde{Y'}$. The K\"{u}nneth formula and
covering space theory lead to the following presentation of
$\pi_2(X')$, viewed as a {\em left}\/ $\Z\pi'$-module:
\begin{equation*}
(\Z\pi')^2
\xrightarrow{
\bigl(\begin{smallmatrix}
1-x_4&1-x_3&0&x_2-1&0&1-x_1&0\\
1-x_5&0&1-x_3&0&x_2-1&0&1-x_1
\end{smallmatrix}\bigr)
}
(\Z\pi')^7\to \pi_{2}(X')\to 0.
\end{equation*}
Notice that $V_6(\pi_2(X'))=\{t\in (\C^*)^5\mid t_1=t_2=t_3=1\}$ is
a $2$-dimensional subtorus.  That $M=\pi_2(X')$ is not nilpotent
can be seen directly, as follows.  Let
$\widetilde{M}=\pi_2(X')\otimes_{\Z\pi'} \Z\Z^5$, and let
$\gr \widetilde{M}$ be the associated graded module (with
respect to the $I$-adic filtration).  From the above presentation,
we may easily compute its Hilbert series:
$\Hilb(\gr \widetilde{M},t)=\frac{7-2t}{(1-t)^5}$.
Since this series is not a polynomial, $ \widetilde{M}$
is not nilpotent, and so $M$ isn't, either.
\end{example}

\section{Graphic arrangements}
\label{sec:graphics}

In this section, we apply our methods to graphic arrangements.
We start by giving a graph-theoretic characterization of hypersolvable
arrangements within this class. We then show how
to get information on
$\pi_2$ of the complement directly from the graph,
in the case of arrangements associated to graphs without
$3$-cycles.

\subsection{Graphs and arrangements}
\label{subsec:grarr}
Let $G=(\V,\E)$ be a non-empty subgraph of
the complete graph on a finite set of
vertices $\V$. Assume that there are no isolated
vertices in the graph, so that the set of edges $\E$ determines $G$.
All graphs considered in this section will be of this type.

Let $\V=\{1,\dots,m\}$. The {\em graphic arrangement} associated
to $G=(\V,\E)$ is the arrangement in $\C^m$ given by
$\A_{G}=\{\ker(z_i-z_j) \mid \{i,j\}\in \E\}$, see \cite{OT}.
For each edge $e=\{i,j\}$, we will denote by $H_e:=\ker(z_i-z_j)$
the corresponding hyperplane of $\A_G$.

Clearly, an arrangement is graphic if and only if it is a
sub-arrangement of a braid arrangement.
For example, if $G$ is the complete graph on $m$ vertices, then
$\A_{G}=\B_{m}$, the braid arrangement in $\C^{m}$.
If $G$ is a diagram of type $\operatorname{A}_{m}$, then
$\A_{G}$ is a Boolean arrangement.
If $G$ is an $m$-cycle, then $\A_{G}$
is the cone of a generic arrangement.

Many of the usual invariants associated to $\A_G$ can be computed
directly from $G$. For example,
$P_{\A_G}(T)=(-T)^m\chi_{G}(-T^{-1})$, where $\chi_{G}(T)$ is the
chromatic polynomial of $G$, see~\cite{OT}.
Also, an \textbf{nbc}-basis for $\OS(\A_G)$ corresponds to an
\textbf{nbc}-basis for $G$, as follows.
Fix an ordering on the edges, $\E_G=\{e_1<\cdots <e_n\}$,
and denote by $\a_i$ the defining equation of $H_{e_i}$.
Then, $\{\a_{i_1},\dots, \a_{i_r}\}$ is minimally dependent
if and only if $\{e_{i_1},\dots, e_{i_r}\}$ is an $r$-cycle of $G$.
Deleting the highest edge from this cycle yields a broken circuit.
The resulting \textbf{nbc}-basis for $\OS(\A_G)$ is given by
\begin{equation}
\label{eq:nbc}
\{e_K \mid \text{$K$ is a subgraph of $G$ which does not contain
any broken circuit of $G$}\},
\end{equation}
where $e_K:=e_{i_1}\cdots e_{i_s}\in \bigwedge^s(e_1,\dots,e_n)$,
if $\E_K=\{e_{i_1},\dots , e_{i_s}\}$.

\subsection{Supersolvable and hypersolvable graphs}
\label{subsec:sshsgr}
The following results, due to Stanley~\cite{St} and
Fulkerson and Gross~\cite{FG}, tell us how to (easily) recognize
supersolvable arrangements within the class of graphic arrangements.

\begin{theorem}
\label{thm:stfg}
Let $\A_G$ be a graphic arrangement. Then:
\newline
\mbox{\quad}\textup{(Stanley~\cite{St})}
  $\A_G$ is supersolvable if and only if the graph
$G$ is supersolvable, i.e., it has a
{\em supersolvable composition series} of
induced subgraphs, $\emptyset = G_0\subset
G_1\subset \cdots\subset G_{\ell}=G$,
such that:
\begin{alphenum}
\item for each $1\le i\le \ell$, there is a single vertex
in $G_i\setminus G_{i-1}$, say, $v_i$;
\item \label{fg}
the subgraph of $G_i$ induced by $v_i$ and its neighbors in $G_i$ is
complete.
\end{alphenum}
\quad\textup{(Fulkerson and Gross~\cite{FG})}
$G$ has a supersolvable composition
series if and only if every cycle in $G$ of length greater than $3$
has a chord.
\end{theorem}

As a simple example, consider the two graphic arrangements
given by the graphs in Figure~\ref{fig:grapharr}.
Neither graph has a $3$-cycle; each graph has
$4$-cycles, with no chords.
Hence, the two arrangements are not supersolvable.

Stanley's supersolvable test from Theorem~\ref{thm:stfg}
has a hypersolvable analogue. To state it, we first need a definition.

\begin{definition}
\label{def:solvgraph}
A pair of graphs, $(G,K)$, is called a {\em solvable extension} if
$K$ is a subgraph of $G$, with $\emptyset \ne \E_{K}\subsetneq \E_{G}$,
and:
\begin{enumerate}
\item\label{c1}
There is no $3$-cycle in $G$ having two edges from $\E_K$
and one edge from $\E_G\setminus \E_K$.
\item \label{c2}
Either $\E_G\setminus \E_K=\{e\}$, and both endpoints of $e$
are not in $\V_K$, or there exist distinct vertices,
$\{v_1,\dots ,v_k,v\}\subset \V_G$,
with $\{v_1,\dots ,v_k\}\subset \V_K$, such that:
\begin{enumerate}
\item \label{c21}
$K$ contains the complete graph
on $\{v_1,\dots ,v_k\}$, and
\item \label{c22}
$\E_G\setminus \E_K=\{\{v,v_s\}\mid 1\le s \le k\}$.
\end{enumerate}
\end{enumerate}
\end{definition}

\begin{lemma}
\label{lem:solvgraph}
An extension of graphs, $(G,K)$, is solvable if and only if
the corresponding extension of graphic arrangements, $(\A_G,\A_K)$,
is solvable.
\end{lemma}

\begin{proof}
Let $e_1=\{i_1,j_1\}$, $e_2=\{i_2,j_2\}$, $e_3=\{i_3,j_3\}$ be three
distinct edges of $G$. Notice that the corresponding defining
equations, $\{z_{i_r}-z_{j_r}\mid 1\le r\le 3\}$, viewed as
points in $\P(\C^{m*})$, are collinear if and only if
$\{e_r\mid 1\le r\le 3\}$ are the edges of a $3$-cycle.
Using this remark, it is a straightforward exercise to translate
conditions \eqref{I}--\eqref{III} from \S\ref{subsec:hyparr} into
conditions \eqref{c1} and \eqref{c2} from
Definition~\ref{def:solvgraph}.
\end{proof}

\begin{definition}
\label{def:hypergraph}
A graph $G$ is called {\em hypersolvable} if it has a
{\em hypersolvable composition series}, i.e., a chain of subgraphs,
$G_1\subset \cdots \subset G_i\subset G_{i+1} \subset \cdots \subset
G_{\ell}=G$,
such that $G_1$ has a single edge, and $(G_{i+1},G_{i})$ is a solvable
extension, for $i=1,\dots,\ell-1$.
\end{definition}

The class of hypersolvable graphs contains the supersolvable
(or chordal) graphs described in Theorem~\ref{thm:stfg}, and many
others.
For example, the graphs from Figure~\ref{fig:grapharr} are both
hypersolvable,
with composition series $G_i=\{e_1,\dots ,e_i\}$, $1\le i\le 9$, but
not
supersolvable.

\begin{prop}
\label{prop:hyptest}
A graph $G$ is hypersolvable if and only if the graphic arrangement
$\A_G$ is hypersolvable.
  \end{prop}

\begin{proof}
Clearly, $G_1\subset \cdots \subset G_{\ell}$
is a composition series for $G$ if and only if
$\A_{G_1}\subset \cdots \subset \A_{G_{\ell}}$
is a composition series for $\A_{G}$.
\end{proof}

\subsection{Graphs with no $3$-cycles}
\label{subsec:notriangle}
We now analyze in more detail a very simple example:
hypersolvable arrangements
coming from graphs without $3$-cycles, and their second homotopy
group.

\begin{prop}
\label{prop:notriangle}
Let $G$ be a graph with no $3$-cycles, and with edges
$\{e_1, \dots , e_n\}$. Then:
\begin{enumerate}
\item \label{q1}
The graph $G$ is hypersolvable, with composition
series $G_i=\{e_1, \dots , e_i \}$, $1 \le i \le n$.
\item \label{q2}
The arrangement $\A=\A_{G}$ is hypersolvable,
with length $n$ and exponents $\{1,\dots,1\}$.
\item \label{q3}
$\QOS(\A)=\bigwedge^*(e_1,\dots , e_n)$.
\item \label{q4}
$\pi_1(X(\A))=\Z^n$.
\end{enumerate}
\end{prop}

\begin{proof}
There are no collinearity relations among the defining equations
of $\A$, since $G$ has no $3$-cycles.
Part~\eqref{q1} then follows from Definitions
\ref{def:solvgraph} and \ref{def:hypergraph},
Part~\eqref{q2} from Proposition~\ref{prop:hyptest} and \eqref{q1}, and
Part~\eqref{q3} from the definition of the quadratic OS-algebra.

Now set $m=\#\{\text{vertices of $G$}\}$. If $m\le 3$, then
Part~\eqref{q4} is trivially verified. If $m>3$, we may take a generic
$3$-plane $P$ in $\C^m$ with the property that
$\pi_1(X(\A))=\pi_1(X(\A\cap P))$
(by \cite{HL}), and such that $\A$ and $\A\cap P$ have the same
collinearity relations. The decone of the arrangement $\A\cap P$
is thus generic, and so $\pi_1(X(\A\cap P))=\Z^n$ (by \cite{Ha}).
\end{proof}

\begin{theorem}
\label{cor:graphpres}
Let $G$ be a graph with no $3$-cycles. Let $\Sq$ be the set of 
$4$-cycles
of $G$. Set $\A=\A_G$ and $X=X(\A)$. Then:
\begin{itemize}
\item
If $\Sq=\emptyset$, then $\pi_2(X)=0$.
\item
If $\Sq\ne\emptyset$, then $\pi_2(X)$ is non-trivial, with
group of coinvariants isomorphic to
$\Z[\Sq]$, the free abelian group generated by $\Sq$.
\end{itemize}
\end{theorem}

\begin{proof}
Set $n=|\A|$, $\pi=\pi_1(X)$ and $p=p(X)$.
Recall from \S\ref{subsec:osalg} the construction of the
OS-algebra of $\A$,
together with the graphic counterpart from \S\ref{subsec:grarr}.
 From Proposition~\ref{prop:notriangle}, we know that $\A$
hypersolvable, $\QOS(\A)=\bigwedge^*(e_1,\dots , e_n)$,
and $\pi=\Z^n$.

If $\Sq=\emptyset$, then $p\geq 3$, by 
Corollary~\ref{cor:pcomb}\eqref{i1},
and thus $\pi_2(X)=0$, by 
Theorem~\ref{thm:hspipres}\eqref{p1}-\eqref{p2}.

If $\Sq\ne\emptyset$, the same argument shows that $p=2$.
 From Theorem~\ref{thm:hspipres}\eqref{p3} and \eqref{eq:dualpr},
we deduce that $\pi_2(X)_{\pi}$ is isomorphic to
$\ker(\pi^3_{\A})$, which is generated by
$\{\partial e_{S} \mid S\in \Sq\}$. A quick inspection
of the construction of the {\bf{nbc}}-basis in degree $3$
reveals that $\ker(\pi^3_{\A})$ is free abelian, with rank equal to
the number of broken $4$-circuits.
Now, in the case of graphic arrangements,
the associated broken circuit uniquely determines a cycle in the graph.
Hence, the map $\partial \colon \Z[\Sq]\to \ker(\pi^3_{\A})$
is actually an isomorphism.
\end{proof}

\begin{example}
\label{ex:pi2ex}
Consider the hypersolvable graphs without $3$-cycles, $G_1$ and $G_2$, 
in
Figure~\ref{fig:grapharr},
and let $\A_1$ and $\A_2$ be the corresponding graphic arrangements
in $\C^7$, with complements $X_1$ and $X_2$. Both complements have
$\pi_1=\Z^{9}$ and
$b_2=\binom{9}{2}=36$. The $4$-cycles of $G_1$ are
$\{ e_1, e_2, e_9, e_7 \}$ and $\{ e_5, e_6, e_7, e_8 \}$,
while those of $G_2$ are $\{ e_1, e_2, e_9, e_7 \}$,
$\{ e_1, e_8, e_6, e_7 \}$ and $\{ e_2, e_9, e_6, e_8 \}$. It follows 
that
the spaces $X_1$ and $X_2$ have distinct homotopy $2$-types, since the 
ranks
of the coinvariants of the $\Z\Z^9$-modules $\pi_2(X_1)$ and
$\pi_2(X_2)$ are different.

\begin{figure}
\setlength{\unitlength}{0.6cm}
\begin{picture}(10,5.5)(2.5,-3.5)
\xygraph{!{0;<8mm,0mm>:<0mm,9.5mm>::}
[]*-{\bullet}
(-^{e_1}[rr]*-{\bullet}
(-^{e_2}[dr]*-{\bullet}
(-^{e_3}[ddl]*-{\bullet}
)
)
,-^{e_8}[ddd]*-{\bullet}
,-_{e_7}[dl]*-{\bullet}
(-_{e_6}[d]*-{\bullet}
(-_{e_5}[dr]*-{\bullet}
(-_{e_4}[rr]*-{\bullet})
)
,-^{e_{9}}[rrrr]*-{\bullet}
)
)
}
\qquad
\xygraph{!{0;<8mm,0mm>:<0mm,9.5mm>::}
[]*-{\bullet}
(-^{e_1}[rr]*-{\bullet}
(-^{e_2}[dr]*-{\bullet}
(-^{e_3}[ddl]*-{\bullet}
)
,-^{e_8\qquad}[ddlll]*-{\bullet}
)
,-_{e_7}[dl]*-{\bullet}
(-_{e_6}[d]*-{\bullet}
(-_{e_5}[dr]*-{\bullet}
(-_{e_4}[rr]*-{\bullet})
)
,-^{\qquad e_{9}}[rrrr]*-{\bullet}
)
)
}
\end{picture}
\caption{\textsf{The graphs $G_1$ and $G_2$}}
\label{fig:grapharr}
\end{figure}

\end{example}
\bigskip

\end{document}